\documentclass[11pt]{article}
\usepackage{latexsym}
\usepackage{amsfonts,amssymb,amsmath}

\usepackage[nodayofweek]{datetime}

\newtheorem{thm}{Theorem}[section]
\newtheorem{cor}[thm]{Corollary}
\newtheorem{lem}[thm]{Lemma}
\newtheorem{pro}[thm]{Proposition}
\newtheorem{defn}[thm]{Definition}

\newcommand{\ov }{\overline }
\title{ Semigroups and one-way functions }    

\author{ J.C.\ Birget  } 

\date{\footnotesize{28 Nov.\ 2015}}
\begin{document}
\maketitle

\begin{center}
{\it To Stuart Margolis on his 60th birthday.} 
\end{center}

\bigskip

\begin{abstract}
We study the complexity classes {\sf P} and {\sf NP} through a semigroup
{\sf fP} (``polynomial-time functions''), consisting of all polynomially 
balanced polynomial-time computable partial functions. 
The semigroup {\sf fP} is non-regular iff ${\sf P} \neq {\sf NP}$.
The one-way functions considered here are based on worst-case complexity
(they are not cryptographic); they are exactly  the non-regular elements 
of {\sf fP}.
We prove various properties of {\sf fP}, e.g., that it is finitely 
generated.  We define reductions with respect to which certain universal 
one-way functions are {\sf fP}-complete. 
\end{abstract}


\smallskip


% Section 1

\section{Introduction}

The goal of this work is to study the complexity classes {\sf P} and 
{\sf NP} via functions, and via semigroups of functions, rather than 
just as sets of languages.
This approach is intuitive (in particular, because of the immediate 
connection with certain one-way functions), and quickly leads to results. 
It is not clear whether this will contribute to a solution of the famous 
{\sf P}-vs.-{\sf NP} problem, but the semigroups considered here, as well 
as the ``inversive reductions'' and the accompanying completeness results 
for one-way functions, are interesting in their own right. 
 
The starting point is a certain kind of {\em one-way functions}, and the 
well-known fact that one-way functions of this kind exist iff 
${\sf P} \neq {\sf NP}$. 

First, some notation: 
We fix an alphabet $A$, which will be $\{0,1\}$ unless another alphabet
is explicitly mentioned. The set of all strings over $A$, including the
empty string, is denoted by $A^*$.
For a partial function $f$, the domain (i.e., the inputs $x$ for which
$f(x)$ is defined) is denoted by ${\sf Dom}(f)$. The image set of $f$ is
denoted by ${\sf Im}(f)$ or by $f(A^*)$ or $f({\sf Dom}(f))$. As a rule we
will use partial functions, even when the word ``partial'' is omitted; we
say {\it total function} for functions whose domain is $A^*$.  
As usual,
{\sf P} and {\sf NP} are the class of languages accepted by deterministic,
respectively nondeterministic, polynomial-time Turing machines
\cite{DuKo,Papadim}.

\medskip

\noindent {\bf Definition scheme:} A function $f: A^* \to A^*$ is called 
{\bf ``one-way''} \ iff \ from $x$ and a description of $f$ it is 
{\sl ``easy''} to compute $f(x)$, but from $f$ and $y \in {\sf Im}(f)$ it 
is {\sl ``difficult''} to find any $x \in A^*$ such that $f(x) = y$. 

\medskip

This is an old idea going back at least to W.S.\ Jevons in 1873, who also 
compared the difficulties of multiplication and factorization of integers 
(as pointed out in \cite{Handbook}). 
The concept became well-known after the work of Diffie and Hellman
\cite{DiffHell}.  
Levin's paper \cite{LevinTale} discusses some deeper connections of one-way 
functions.
The definition scheme can be turned into precise definitions, in many
(non-equivalent) ways, by defining ``easy'' and ``difficult'' (and, if needed, 
``description'' of $f$).

\begin{defn}
A partial function $f: A^* \to A^*$ is {\em polynomially balanced} iff 
there exists polynomials $p, q$ such that for all $x \in {\sf Dom}(f)$:
 \, $|f(x)| \leq p(|x|)$ and $|x| \leq q(|f(x)|)$.
\end{defn}
We call the polynomial $q$ above an {\em input balance} function of $f$.
The word ``honest'' is often used in the literature for polynomially 
balanced.

We introduce the following set of ``easy'' functions.

\begin{defn} {\bf (the semigroup {\sf fP}).} 
 \ {\sf fP} is the set of partial functions $f: A^* \to A^*$ that are 
polynomially balanced, and such that 
 \ $x \in {\sf Dom}(f) \longmapsto f(x)$ \ is computable by a 
deterministic polynomial-time Turing machine. 
(It follows from the second condition that ${\sf Dom}(f)$ is in {\sf P}.)
\end{defn}
As a rule, a machine that computes a partial function $f$ can always also
be used as an acceptor of ${\sf Dom}(f)$.  

When $A$ is an arbitrary alphabet (as opposed to $\{0,1\}$) we write 
${\sf fP}_A$ or ${\sf fP}_{|A|}$.
The complexity class {\sf fP} is different from the complexity class 
{\sf FP}, considered in the literature \cite{Papadim}; {\sf FP} is a set 
of relations (viewed as search problems) whereas {\sf fP} is a set of 
partial functions. 
It is easy to see that ${\sf fP}_{|A|}$ is closed under composition, so it 
is a semigroup.

\begin{defn} {\bf (worst-case deterministic one-way function).}  
A partial function $f: A^* \to A^*$ is {\em one-way} iff $f \in {\sf fP}$, 
but there exists {\em no} deterministic polynomial-time algorithm which, on 
every input $y \in {\sf Im}(f)$, 
outputs some $x \in A^*$ such that \ $f(x) = y$.
There is no requirement when $y \not\in {\sf Im}(f)$.
\end{defn} 
This kind of one-way functions is defined in terms of worst-case complexity, 
hence it is not ``cryptographic''. However, it is relevant to the 
{\sf P}-vs.-{\sf NP} problem because of the following fact (see e.g., 
\cite{HemaOgi} p.\ 33 for a proof and history).  

\begin{pro} {\bf (folklore).} 
 \ One-way functions (in the worst-case sense) exist iff 
 \ ${\sf P} \neq {\sf NP}$. \ \ \ $\Box$
\end{pro}
The concept of an {\it inverse} is central to one-way functions.
The following lemma is straightforward to prove.

{\sf Notation:} 
For a partial function $f$ and a set $S$, the restriction of $f$ to $S$ is 
denoted by $f|_S$; for the restriction of the identity map to $S$ we simply 
write ${\sf id}_S$.

\begin{lem} \label{concept} {\bf (concept of inverse).} 

\noindent
For partial functions $f, f': A^* \to A^*$ the following are equivalent.

\smallskip

\noindent $\bullet$ \ For all $y \in {\sf Im}(f)$, \ $f'(y)$ is defined
and $f(f'(y)) = y$.  

\smallskip

\noindent $\bullet$ \ $f \circ f'|_{{\sf Im}(f)} = {\sf id}_{{\sf Im}(f)}$ .

\smallskip

\noindent $\bullet$ \ $f \circ f' \circ f = f$. \ \ 

\smallskip

\noindent These properties imply ${\sf Im}(f) \subseteq {\sf Dom}(f')$.  
 \ \ \ \ \  $\Box$
\end{lem}

\begin{defn} \ A function $f'$ such that $f \circ f' \circ f = f$ is called 
an {\em inverse of} $f$.
\end{defn}
The following recipe gives more intuition about inverses.

\medskip

\noindent Pseudo-algorithm: {\it How any inverse $f'$ of a given $f$ is 
made.}

\smallskip

\noindent (1) \ Choose ${\sf Dom}(f')$ such that 
${\sf Im}(f) \subseteq {\sf Dom}(f')$. 

\smallskip

\noindent (2) \ For every $y \in {\sf Im}(f)$, choose $f'(y)$ to be any
$x \in f^{-1}(y)$.

\smallskip

\noindent (3) \ For every $y \in {\sf Dom}(f') - {\sf Im}(f)$, choose 
$f'(y)$ arbitrarily in $A^*$.  

\bigskip

\noindent {\bf Remark.} When $f'$ is an inverse of $f$, the restriction 
$f'|_{ {\sf Im}(f)}: y \in {\sf Im}(f) \mapsto f'(y)$ is the {\em choice
function} corresponding to $f'$. For set theory in general, the existence 
of choice functions (and the existence of inverses) for every partial 
function is equivalent to the Axiom of Choice. The existence of one-way 
functions in our sense amounts to the non-existence of certain inverses; 
the existence of one-way functions is thus equivalent to the {\em non-validity 
of the Axiom of Choice} in the (highly restricted) context of deterministic 
polynomial time-complexity. 

From the definition of polynomially balanced we can see now:
{\em If $f$ is polynomially balanced then so is every choice function 
corresponding to any inverse of $f$. }

\begin{defn} \ Let $S$ be a semigroup. An element $x \in S$ is called
{\em regular} iff there exists $x' \in S$ such that $x x' x = x$. 
In that case, $x'$ is called an {\em inverse} of $x$. 
A semigroup $S$ is called regular iff every element of $S$ is regular.

Let $S$ be a monoid with identity {\bf 1}. Then $x'$ is a {\em left}- (or 
{\em right}-) inverse of $x$ iff $x'x = {\bf 1}$ (respectively 
$x x' = {\bf 1}$). If $x'x = x x' = {\bf 1}$ then $x'$ is a two-sided inverse
or group-inverse. {\rm (See \cite{CliffPres,Grillet}.)}
\end{defn}

\noindent The following summarizes what we have seen, and gives the
initial motivation for studying the class {\sf NP} via certain functions
and semigroups. 

\begin{pro} \  
The monoid {\sf fP} is {\em regular} iff one-way functions do {\sl not} 
exist (iff \ ${\sf P} = {\sf NP}$).  \ \ \ \ \ $\Box$
\end{pro} 
Some properties of the image set of functions in {\sf fP}:

\begin{pro}  \label{ImfP} \!\!\! .

\noindent (1) \ For every $f \in {\sf fP}$, \, ${\sf Im}(f) \in {\sf NP}$.

\noindent (2) \ If $f \in {\sf fP}$ and $f$ is regular then
${\sf Im}(f) \in {\sf P}$.

\noindent (3) \ For every language $L \in {\sf NP}$ there exists
$f_L \in {\sf fP}$ such that $L = {\sf Im}(f_L)$. 

 \ Moreover, the set of functions $\{f_L \in {\sf fP} : L \in {\sf NP}\}$ 
can be chosen so that $f_L$ is regular iff \ $L \in {\sf P}$.
The map $L \in {\sf NP} \mapsto f_L \in {\sf fP}$ is an embedding of
{\sf NP} (as a set) into {\sf fP}, such that {\sf P} (and only {\sf P}) is 
mapped into the regular elements of {\sf fP}. 
\end{pro}
\noindent {\bf Proof.} (1) is obvious (polynomial balance is needed).

\noindent (2)  Let $f' \in {\sf fP}$ be an inverse of $f$.
If $y \in {\sf Im}(f)$ then $ff'(y) = y$.
If $y \not\in {\sf Im}(f)$, then either $f'(y)$ is not defined, or
$ff'(y) \in {\sf Im}(f)$, hence $ff'(y) \neq y$. Thus, $y \in {\sf Im}(f)$ 
iff $ff'(y) = y$.
When $f, f' \in {\sf fP}$ then on input $y$ the properties $ff'(y) = y$,
$y \not\in {\sf Dom}(f')$, $ff'(y) \neq y$, can be decided deterministically
in polynomial time.

\noindent (3) Let $M$ be a nondeterministic Turing machine accepting $L$, 
such that all computations of $M$ are polynomially bounded, and do not halt
before the whole input has been read. We can assume that $M$ has binary 
nondeterminism, i.e., in each transition it has at most two nondeterministic 
choices.  Define 

\smallskip

 \ \ \ \ \ \   
 $f_L(x,s) = x$ \ \ iff \ \ $M$, with {\it choice sequence} $s$, 
 accepts $x$; 

 \ \ \ \ \ \    
 $f_L(x,s)$ is undefined otherwise. 

\smallskip 

\noindent Choice sequences are also called guessing sequences, or advice 
sequences.
Then $L = {\sf Im}(f_L)$ and $f_L \in {\sf fP}$; balancing comes from the 
fact that $s$ has polynomially bounded length.

We saw that if $f_L$ is regular then ${\sf Im}(f_L) = L \in {\sf P}$.
Moreover, if $L \in {\sf P}$ then the Turing machine $M$ can be chosen to
be deterministic, and the choice sequence $s$ is the all-0 word;
so in that case, $f_L$ is not one-way.
  \ \ \ $\Box$

\begin{cor}  \ The {\em image transformation} {\sf Im}: 
$f \mapsto {\sf Im}(f)$ maps {\sf fP} onto {\sf NP}, and maps the set of 
regular elements of {\sf fP} onto {\sf P}. 
Moreover, {\sf NP} is embedded into {\sf fP} (by the transformation 
$L \mapsto f_L$ above), and {\sf NP} is a {\em retract} of {\sf fP} (by the 
transformations $L \mapsto f_L$ and {\sf Im}).
 \ \ \   \ \ \ $\Box$
\end{cor}
The following suggests that ${\sf Im}(f) \in {\sf P}$ is not equivalent 
to the regularity of $f$.

\begin{thm} 
 \ If $\Pi_2^{\sf P} \neq \Sigma_2^{\sf P}$ then there exist surjective 
one-way functions. {\rm (Proved in \cite{JCBcircsizeinv}.)} 
\end{thm}

In the next sections we show various properties of {\sf fP} and of 
closely related semigroups. The fact that these semigroups have 
interesting properties, and that the proofs are not difficult, gives 
a second motivation for studying these semigroups.

\section{The Green relations of polynomial-time function semigroups}

We give a few properties of the Green relations $\leq_{\cal R}, 
 \ \leq_{\cal L}, \ \leq_{\cal J}, \ \equiv_{\cal D}$ (see
\cite{CliffPres,Grillet}).  First some notation.

Let $F: X \to Y$ be a partial function; then $F^{-1}: Y \to X$ denotes 
the {\em inverse relation} of $F$, i.e., for all 
$(x,y) \in X \times Y$: \ $x \in F^{-1}(y)$ iff $F(x) = y$. By
${\sf mod}F$ we denote the partition on ${\sf Dom}(F)$, defined by
\ $x_1$ ${\sf mod}F$ $x_2$ iff $F(x_1) = F(x_2)$.
The set of blocks (equivalence classes) of ${\sf mod}F$ is 
\ $\{ F^{-1}F(x) : x \in {\sf Dom}(F)\}$.
For two partial functions $F, C: X \to Y$ we write
${\sf mod}C \leq {\sf mod}F$ (the partition of $C$ is {\em coarser} 
than the partition of $F$, or the partition of $F$ {\em refines} the 
partition of $C$) iff \ ${\sf Dom}(C) \subseteq {\sf Dom}(F)$ \ and
 \ $F^{-1}F(x) \subseteq C^{-1}C(x)$ \ for all $x \in {\sf Dom}(C)$;
equivalently, every ${\sf mod}C$-class is a union of ${\sf mod}F$-classes.

\begin{pro}{\rm (regular ${\cal L}$- and $\cal R$-orders).}
\label{regLR}   
If $f, r \in {\sf fP}$ and $r$ is {\em regular} with an inverse
$r' \in {\sf fP}$ then:

\smallskip

\noindent $\bullet$ \ $f \leq_{\cal R} r$ \ \ iff \ \ $f = r r' f$
 \ \ iff \ \ ${\sf Im}(f) \subseteq {\sf Im}(r)$.

\smallskip

\noindent $\bullet$ \ $f \leq_{\cal L} r$ \ \ iff \ \ $f = fr'r$  
 \ \ iff \ \ ${\sf mod}f \leq {\sf mod}r$. 
\end{pro}
{\bf Proof.} [$\cal R$-order]:  \ 
$f \leq_{\cal R} r$ \ iff \ for some $u \in {\sf fP}: f = ru$. Then 
$f = r r' r u = r r' f$. 
Also, it is straightforward that $f = ru$ implies 
 ${\sf Im}(f) \subseteq {\sf Im}(r)$.

Conversely, if ${\sf Im}(f) \subseteq {\sf Im}(r)$ then 
 \ ${\sf id}_{{\sf Im}(f)} \ = \ $
${\sf id}_{{\sf Im}(r)} \circ {\sf id}_{{\sf Im}(f)} \ = \ $
$r \circ r'|_{{\sf Im}(r)} \circ {\sf id}_{{\sf Im}(f)}$. Hence,
$f \ = \ {\sf id}_{{\sf Im}(f)} \circ f \ = \ $
$r \circ r'|_{{\sf Im}(r)} \circ {\sf id}_{{\sf Im}(f)} \circ f \ = \ $
$r \circ r'|_{{\sf Im}(r)} \circ f \ \leq_{\cal R} r$. 
 
\medskip

[$\cal L$-order]: \
$f \leq_{\cal L} r$ \ iff \ for some $v \in {\sf fP}: f = vr$. Then 
$f = v r r' r = f r' r$. 
And it is straightforward that $f = vr$ implies
 ${\sf mod}f \leq {\sf mod}r$.

Conversely, if ${\sf mod}f \leq {\sf mod}r$ then for all 
$x \in {\sf Dom}(f)$, \ $r^{-1}r(x) \subseteq f^{-1}f(x)$.  And for 
every $x \in {\sf Dom}(f)$, \ $\{f(x)\} = f \circ f^{-1} \circ f(x)$.
Moreover,
$f \circ r^{-1} \circ r(x) \subseteq f \circ f^{-1} \circ f(x) = \{f(x)\}$, 
and since $r^{-1} \circ r(x) \neq \varnothing$, it follows that 
$f \circ r^{-1} \circ r(x) = \{f(x)\}$.  So, $f = f \circ r^{-1} \circ r$. 
Moreover, $f \circ r' \circ r(x) \in f \circ r^{-1} \circ r(x) = \{f(x)\}$,
hence $f \circ r' \circ r(x) = f(x)$.  Hence, $f = f r' r \leq_{\cal L} r$.  
 \ \ \ $\Box$
 
\bigskip

\noindent The ${\cal D}$-relation between elements of {\sf fP} with infinite
image sets seems difficult, even in the case of regular elements.
A first question (inspired from the Thompson-Higman monoids 
\cite{JCBmonThH}): Are all regular elements of {\sf fP} with infinite image 
in the same ${\cal D}$-class, i.e., the ${\cal D}$-class of 
${\sf id}_{A^*}$?

\begin{pro} \label{D_inj}  Let $f \in {\sf fP}$ be regular. 
Then $f \equiv_{\cal D} {\sf id}_{A^*}$ iff there exists $g \in {\sf fP}$ 
such that $g$ is injective, total, and regular, and such that 
${\sf Im}(f) = {\sf Im}(g)$.
\end{pro}
{\bf Proof.} Assume $f \equiv_{\cal D} {\sf id}_{A^*}$. By Prop.\ \ref{regLR},
$f \equiv_{\cal R} {\sf id}_L$, where $L = {\sf Im}(f)$.
And ${\sf id}_{A^*} \equiv_{\cal D} {\sf id}_L$ iff there exists 
$g \in {\sf fP}$ such that 
${\sf id}_{A^*} \equiv_{\cal L} g \equiv_{\cal R} {\sf id}_L$.
Hence, ${\sf id}_{A^*} = g' \, g$ for some $g' \in {\sf fP}$; this equality
implies that $g$ is total and injective. The existence of $g' \in {\sf fP}$
implies that $g$ is regular. Since $g \equiv_{\cal R} {\sf id}_L$, 
${\sf Im}(g) = L$. Hence, $g$ has the required properties. 

To prove the converse we will use the following:

\medskip

\noindent {\sf Claim.} \ For every $g \in {\sf fP}$ we have:
 \ $g$ is injective, total, and regular \ iff 
 \ $(\exists g' \in {\sf fP}) \, g' g = {\sf id}_{A^*}$. 

\noindent Proof of the Claim. 
The right-to-left implication is straightforward. In the other direction, 
if $g$ is regular then there exists $g' \in {\sf fP}$ such that $gg'g =g$. 
And if $g$ is total and injective, there exists a partial function $h$ such 
that $hg = {\sf id}_{A^*}$. Now $gg'g =g$ implies $hgg'g = hg$, hence by 
using $hg = {\sf id}_{A^*}$ we obtain: $g' g = {\sf id}_{A^*}$. 
This proves the Claim.

\medskip
 
For the converse of the Proposition, assume there exists $g \in {\sf fP}$
with the required properties. 
If such a $g$ exists, then $f \equiv_{\cal R} g$, by Prop.\ \ref{regLR}. 
Moreover, $g \equiv_{\cal L} {\sf id}_{A^*}$; this follows from 
$g' g = {\sf id}_{A^*}$, which holds by the Claim. Hence
$f \equiv_{\cal R} g \equiv_{\cal L} {\sf id}_{A^*}$.
 \ \ \ $\Box$

\bigskip

\noindent However, it is an open problem whether every infinite language $L$ 
in {\sf P} is the image of an injective, total, polynomial-time computable 
function $g$ (and whether $g$ can be taken to be regular or one-way).
Also, not much is known about which infinite languages in {\sf P} can be 
mapped onto each other by maps in {\sf fP}.

When ${\sf Im}(f)$ is a right ideal, more can be said.  
By definition, a {\em right ideal} of $A^*$ is a subset $R \subseteq A^*$
such that $R \, A^* = R$ (i.e., $R$ is closed under right-concatenation by 
any string).
Equivalently, a right ideal is a set of the form $R = L \, A^*$, for any set 
$L \subseteq A^*$; in that case we also say that $L$ {\em generates} $R$
as a right ideal.
A {\em prefix code} in $A^*$ is a set $P \subseteq A^*$ such that no word 
in $P$ is a prefix of another word in $P$.  It is not hard to prove that 
for any right ideal $R$ there exists a unique prefix code $P_{_R}$ such that
$R = P_{_{\!\! R}} \, A^*$; in other words, $P_{_R}$ is the minimum generating set 
of $R$, as a right ideal.

A {\em right ideal morphism} is a partial function $h: A^* \to A^*$ such
that for all $x \in {\sf Dom}(h)$ and all $w \in A^*$: \ $h(xw) = h(x) \, w$.
One proves easily that then ${\sf Dom}(h)$ and ${\sf Im}(h)$ are right 
ideals.  

We also consider $A^{\omega}$ (the $\omega$-sequences over $A$, see e.g.\
\cite{PerrinPin}). For a set $L \subseteq A^*$ we define ${\sf ends}(L)$ to 
consist of all elements of $A^{\omega}$ that have a prefix in $L$.  
The Cantor space topology on $A^{\omega}$ uses the sets of the form 
${\sf ends}(L)$ (for $L \subseteq A^*$) as its open sets; here we can 
assume without loss of generality that $L$ is a prefix code or a right 
ideal of $A^*$.

\begin{lem} \label{codeVSrideal} \
If a right ideal $R \subseteq A^*$ belongs to {\sf P} then the corresponding 
prefix code $P$ (such that $R = PA^*$) also belongs to {\sf P}. Conversely, 
if $L$ is in {\sf P} then $LA^*$ is in {\sf P}.
\end{lem}
{\bf Proof.} The first statement follows immediately from the fact 
that $x \in P$ iff $x \in R$ and every strict prefix of $x$ does not belong
to $R$. The converse is straightforward.   \ \ \ $\Box$

\bigskip

\noindent {\bf Notation.} \ Below, $\ov{PA^*}$ \, denotes $A^* - PA^*$ 
(complement).

\begin{pro} \label{Pp_0} \ Let $P \subseteq A^*$ be a prefix code that 
belongs to {\sf P}, and let $p_0 \in P$. Then all {\em regular} elements 
$r \in {\sf fP}$ whose image is \ ${\sf Im}(r) \, = \, L_P \, = \, $
$(P - \{p_0\}) \, A^* \, \cup \ p_0 \, (p_0 A^* \ \cup \ \ov{PA^*})$ \ are 
in the ${\cal D}$-class of \, ${\sf id}_{A^*}$.
We can view $L_P$ as an ``approximation'' of the right ideal $PA^*$ since 

\medskip

\hspace{1in}  $(P - \{p_0\}) \, A^* \ \subset \ L_P \ \subset \ P A^*$. 
 
\end{pro} 
{\bf Proof.} \ Let $L = L_P = {\sf Im}(r)$. By Prop.\ \ref{regLR}, 
$r \equiv_{\cal R} {\sf id}_L$, so it suffices to prove that 
 \ ${\sf id}_L \equiv_{\cal D} {\sf id}_{A^*}$.
We define $\pi, \pi' \in {\sf fP}$ by
\[ \pi(x) = \left\{  
   \begin{array}{ll}
       x     & \mbox{if  \ $x \in (P - \{p_0\}) A^*$, } \\  
       p_0 x & \mbox{otherwise (i.e., if $x \in p_0A^*$ or $x \not\in PA^*$);}
   \end{array}  \right.
\]
\[ \pi'(x) = \left\{ 
   \begin{array}{ll} 
        x   & \mbox{if \ $x \in (P - \{p_0\}) A^*$, } \\  
        z   & \mbox{if \ $x \in p_0A^*$ with $x = p_0z$, } \\  
        {\rm undefined} & \mbox{otherwise (i.e., when $x \not\in PA^*$). }
   \end{array}  \right.
\]
\noindent Then $\pi$ is a total and injective function, and
${\sf Im}(\pi) = L$. Hence, $\pi \equiv_{\cal R} {\sf id}_L$.
Moreover, $\pi' \circ \pi = {\sf id}_{A^*}$, as is easily verified, 
hence $\pi \equiv_{\cal L} {\sf id}_{A^*}$. In summary,
 \ $r \equiv_{\cal R} {\sf id}_L \equiv_{\cal R} \pi \equiv_{\cal L}$
${\sf id}_{A^*}$. 
  \ \ \ \ \ $\Box$

\bigskip

\noindent 
Functions that have right ideals as domain and image are of great interest, 
because of the remarkable properties of the Thompson-Higman 
groups and monoids \cite{McKTh,Th,Scott,Hig74,CFP,BiThomps,BiCoNP} 
and \cite{JCBonewPerm, JCBmonThH, JCBrl}.
Prop.\ \ref{Pp_0} gives an additional motivation for looking at the special
role of right ideals.
This motivates the following.

\begin{defn} \label{RM2P} 
 \ \ \ ${\cal RM}_{_{|A|}}^{\sf P} \ = \ \{ f \in {\sf fP} \ : \ $
$f$ is a right ideal morphism of $A^* \}$.
\end{defn}
When $f$ is a right ideal morphism, ${\sf Dom}(f)$ and ${\sf Im}(f)$ are
right ideals.
${\cal RM}_{_{|A|}}^{\sf P}$ is closed under composition, and
${\cal RM}_2^{\sf P}$ is a submonoid of {\sf fP}.

An interesting submonoid of ${\cal RM}_{_{|A|}}^{\sf P}$ is 
${\cal RM}_{_{|A|}}^{\sf fin}$, consisting of all those
$f \in {\cal RM}_{_{|A|}}^{\sf P}$ for which ${\sf Dom}(f)$ (and hence also
${\sf Im}(f)$) is a {\it finitely generated} right ideal. The monoid
${\cal RM}_{_{|A|}}^{\sf fin}$ is used to define the Thompson-Higman monoid
$M_{|A|,1}$ in \cite{JCBmonThH}. 

\begin{pro} \label{fPinvinRM} 
If an element $f \in {\cal RM}_2^{\sf P}$ has an inverse in {\sf fP} then
$f$ also has an inverse in ${\cal RM}_2^{\sf P}$.
\end{pro}
{\bf Proof.} Let $f_0' \in {\sf fP}$ be an inverse of $f$; we want to 
construct an inverse $f'$ of $f$ that belongs to ${\cal RM}_2^{\sf P}$.
Since $f$ is regular in {\sf fP}, we know from Prop.\ \ref{ImfP} that 
${\sf Im}(f)$ is in {\sf P}.
Hence we can restrict $f_0'$ to ${\sf Im}(f)$, i.e., 
${\sf Dom}(f_0') = {\sf Im}(f)$.
We proceed to define $f'(y)$ for $y \in {\sf Im}(f)$.

First, we compute the shortest prefix $p$ of $y$ that satisfies 
$p \in {\sf Dom}(f_0') = {\sf Im}(f)$. Since  ${\sf Im}(f) \in {\sf P}$, 
this can be done in polynomial time. Now, $y = p \, z$ for some  string $z$. 

Second, we define $f'(y) \ = \ f'_0(p) \ z$, \ where $p$ and $z$ are as 
above.  Thus, $f'$ is a right-ideal morphism.

\smallskip

Let us verify that $f'$ has the claimed properties.
Clearly, $f'$ is polynomial-time computable, and polynomially 
balanced (the latter following from the fact that $f'$ is an inverse of $f$,
which we prove next).
To prove that $f'$ is an inverse of $f$, let $x \in {\sf Dom}(f)$.
Then $f (f' (f(x))) = f (f'(p \, z))$, where $y = f(x) = p \, z$, and $p$
is the shortest prefix of $y$ such that $p \in {\sf Im}(f)$.
Then, $f'(p \, z) = f_0'(p) \, z$, by the definition of $f'$.
Then, since $f$ is a right-ideal morphism, $f(f_0'(p) \, z)$
$= f(f_0'(p)) \, z = p \, z$ (the latter since $f'_0$ is an inverse of $f$,
and since $p \in {\sf Im}(f)$).
Hence, $ff'|_{{\sf Im}(f)} = {\sf id}_{{\sf Im}(f)}$. Thus, by Prop.\ 
\ref{concept}, $f'$ is an inverse of $f$. 
 \ \ \ $\Box$

\bigskip

\noindent {\bf Remark and notation:} 
 \ For $f \in {\cal RM}_{_{|A|}}^{\sf P}$ we saw that ${\sf Dom}(f)$ and 
${\sf Im}(f)$ are right ideals. Let ${\sf domC}(f)$, called the 
{\em domain code}, be the prefix code that generates ${\sf Dom}(f)$ as a 
right ideal. Similarly, let ${\sf imC}(f)$, called the {\em image code}, 
be the prefix code that generates ${\sf Im}(f)$.

In general, ${\sf imC}(f) \subseteq f({\sf domC}(f))$, and it can happen that
${\sf imC}(f) \neq f({\sf domC}(f))$. However the last paragraph of proof of
Prop.\ \ref{fPinvinRM} shows that in any case: 
 \ {\em If $f \in {\cal RM}_{_{|A|}}^{\sf P}$ is regular then $f$ has an 
inverse $f' \in {\cal RM}_{_{|A|}}^{\sf P}$ such that 
${\sf domC}(f') = {\sf imC}(f)$. }

\bigskip

\noindent {\bf Notation:} \ For two words $u,v \in A^*$,
$(v \leftarrow u)$ denotes the right ideal morphism $ux \mapsto vx$
(for all $x \in A^*$). In particular,
$(\varepsilon \leftarrow \varepsilon) = {\sf id}_{A^*}$, where
$\varepsilon$ denotes the empty word.
The morphism $(v \leftarrow u)$ is length-balanced because $|u|, |v|$ are
constants for a given morphism.

\begin{pro} \label{j0} \ For every alphabet $A$, the monoid 
${\cal RM}_{_{|A|}}^{\sf P}$ is ${\cal J}^0$-simple (i.e., the only ideals
are $\{0\}$ and ${\cal RM}_{_{|A|}}^{\sf P}$ itself).
\end{pro}
{\bf Proof.} For any $f \in {\cal RM}_{|A|}^{\sf P}$ that is not the 
empty map, there exist words $x_0, y_0$ such that $f(x_0) = y_0$. Then 
$(\varepsilon \leftarrow \varepsilon)$ $=$ 
$(\varepsilon \leftarrow y_0) \circ f \circ (x_0 \leftarrow \varepsilon)$.
Hence, ${\sf id}_{A^*} \leq_{\cal J} f$ for every non-empty element 
$f \in {\cal RM}_{|A|}^{\sf P}$.  
 \ \ \ $\Box$

\begin{pro} \ \ {\sf fP} is not ${\cal J}^0$-simple, and it has regular 
elements in different non-0 ${\cal J}$-classes.
\end{pro}
{\bf Proof.} The map $\ell: x \in \{0,1\}^* \longmapsto 0^{|x|}$
is in {\sf fP} and it is an idempotent. 

Moreover, $\ell \not\equiv_{\cal J} {\sf id}_{A^*}$.  
Indeed, if there exist functions $\beta, \alpha$ such that for all 
$x \in A^*$, $x = \beta \, \ell \, \alpha(x) = \beta(0^{|\alpha(x)|})$, then 
$|\alpha(x)|$ is different for every $x \in A^*$. But then $\alpha$ is not
polynomially balanced, since $|\alpha(x)|$ would have to range over 
$|A|^{|x|}$ values.  
 \ \ \ $\Box$

\begin{cor} \ \ {\sf fP} and ${\cal RM}_{_{|A|}}^{\sf P}$ are not
isomorphic. \ \ \  \ \ \ $\Box$
\end{cor}

\noindent  As a consequence of Prop. \ref{Pp_0} we have:

\begin{cor} \ Every {\em regular} element $r \in {\cal RM}_2^{\sf P}$ is 
``close'' to an element of {\sf fP} belonging to the ${\cal D}$-class of
${\sf id}_{A^*}$. 
Here, $h_{p_0} \in {\sf fP}$ is  called ``close'' to $r$ \, iff \  
${\sf Im}(r) = PA^*$ for a prefix code $P$, and there exists $p_0 \in P$ such
that:

\smallskip

\noindent $\bullet$ \ 
$(P - \{p_0\}) \, A^* \ \subseteq \ {\sf Im}(h_{p_0}) \ \subseteq PA^*$, 
 \ and 

\smallskip

\noindent $\bullet$ 
 \ $h_{p_0}(x) = r(x)$ whenever $r(x) \in {\sf Im}(h_{p_0})$. 

\end{cor}
{\bf Proof.} Let $P = {\sf domC}(r)$, so $PA^* = {\sf Im}(r)$. 
For every $p_0 \in P$, $r$ is close to ${\sf id}_{L_P} \circ r$, whose image 
set is $L_P \ = \ $ 
$(P - \{p_0\}) \, A^* \ \cup \ p_0 \, (p_0 A^* \cup \ov{PA^*})$, hence
 \ $(P - \{p_0\}) \, A^* \subset L_P \subset PA^*$. 
And ${\sf id}_{L_P} \circ r \equiv_{\cal R} {\sf id}_{L_P}$ since
$L_P \subset PA^*$. 
 \ \ \ $\Box$

\medskip

\noindent Recall the notation $(v \leftarrow u)$ given just before Prop.\
\ref{j0}.

\begin{lem}  \label{RM1LR} \   
In ${\cal RM}_2^{\sf P}$, the $\cal L$-class of ${\sf id}_{A^*}$ is 
$\{(v \leftarrow \varepsilon) \in {\cal RM}_2^{\sf P} : v \in A^*\}$.
This is the set of elements of ${\cal RM}_2^{\sf P}$ that are injective 
and total (i.e., defined for all $x \in A^*$).

The $\cal R$-class of ${\sf id}_{A^*}$ is 
$\{ f \in {\cal RM}_2^{\sf P} : \,  \varepsilon \in {\sf Im}(f)\}$. 
This is the set of elements of ${\cal RM}_2^{\sf P}$ that are surjective 
(i.e., map onto $A^*$).
\end{lem}
{\bf Proof.} If $f \equiv_{\cal L} {\sf id}_{A^*}$ then 
$\varepsilon \in {\sf Dom}(f) = A^*$, so there is $v \in A^*$ such that 
$v = f(\varepsilon)$. Then $f(x) = vx$ for all $x\in A^*$. 
Conversely, if $f(x) = vx$ for all $x\in A^*$ then 
$(\varepsilon \leftarrow v) \circ f = {\sf id}_{A^*}$.

If $f \equiv_{\cal R} {\sf id}_{A^*}$ then ${\sf Im}(f) = A^*$. So
$\varepsilon \in {\sf Im}(f)$. 
Conversely, if $f$ satisfies $\varepsilon \in {\sf Im}(f)$, i.e., 
$\varepsilon = f(x_0)$ for some $x_0 \in A^*$, then 
$f \circ (x_0 \leftarrow \varepsilon) = (\varepsilon \leftarrow \varepsilon)$
$= {\sf id}_{A^*}$.
 \ \ \ $\Box$

\bigskip

\noindent Lemma \ref{RM1LR} shows that the $\cal L$-class of
${\sf id}_{A^*}$ in ${\cal RM}_2^{\sf P}$ is also the $\cal L$-class of
${\sf id}_{A^*}$ in ${\cal RM}_2^{\sf fin}$.

\begin{pro} \ ${\cal RM}_2^{\sf P}$ has trivial group of units, i.e., the 
$\cal D$-class of the identity ${\sf id}_{A^*}$ is $\cal H$-trivial.
\end{pro}
{\bf Proof.} 
If $f \equiv_{\cal H} {\sf id}_{A^*}$ then by Lemma \ref{RM1LR} (for the 
$\cal L$-class of ${\sf id}_{A^*}$), $f(x) = vx$ for all $x$. 
Also by Lemma \ref{RM1LR} (for the $\cal R$-class of ${\sf id}_{A^*}$), 
$f(x_1) = v x_1 = \varepsilon$, for some $x_1$.
This implies $v = \varepsilon$, hence $f = {\sf id}_{A^*}$.
 \ \ \ $\Box$

\bigskip

\noindent  As a consequence of Lemma \ref{RM1LR}, ${\cal RM}_2^{\sf P}$ can 
be injectively mapped (non-homomorphically) into the $\cal R$-class of 
${\sf id}_{A^*} \in {\cal RM}_2^{\sf P}$.  Let us define $f \mapsto \psi_f$ 
by ${\sf Dom}(\psi_f) = \{0\} \, \cup \, 1 \, {\sf Dom}(f)$, and

\smallskip

 \ \ \ \ \ $\psi_f(0) \ = \ \varepsilon$, \ \ and

\smallskip

 \ \ \ \ \ $\psi_f(1x) \ = \ 1 \, f(x)$, \, for all $x \in {\sf Dom}(f)$.

\smallskip

\noindent Then for all $1x \in 1 \, A^*$, 
 \ $\psi_{g \circ f}(1x) = (\psi_g \circ \psi_f)(1x) = 1 \, (f g)(x)$.
So, $f \mapsto \psi_f|_{1A^*}$ is a morphism (where $\psi_f|_{1A^*}$ is
the restriction of $\psi_f$ to $1A^*$). 
But $\psi$ is not a morphism; indeed, since ${\cal RM}_2^{\sf P}$ 
contains non-trivial groups, but the $\cal D$-class of ${\sf id}_{A^*}$ is 
$\cal H$-trivial, there cannot be a homomorphic embedding of 
${\cal RM}_2^{\sf P}$ into the $\cal D$-class of ${\sf id}_{A^*}$.

%%%%%%%%%%%%%%%%%%%%%%%%%%%%%%%%%%%%%%%%%%%%%%%%%%%%%%%%%%%%%%%%%%%%%%%%

\section{\bf Embedding {\sf fP} into ${\cal RM}_2^{\sf P}$ }

\noindent {\bf Transforming any map into a right-ideal morphism:}

\medskip

\noindent The semigroup {\sf fP} uses the alphabet $\{0,1\}$; let $\#$ be 
a new letter.  For $f \in {\sf fP}$ we define
$f_{\#}: \{0, 1, \#\}^* \to \{0, 1, \#\}^*$ \, by letting \,  
${\sf Dom}(f_{\#}) \ = \ {\sf Dom}(f) \, \# \, \{0, 1, \#\}^*$, and

\smallskip

 \ \ \ \ \ \ $f_{\#}(x \# w) \ = \ f(x) \ \# \, w$, 

\smallskip

\noindent for all $x \in {\sf Dom}(f)$ $( \subseteq \{0,1\}^*)$, and all 
$w \in \{0, 1, \#\}^*$. So ${\sf domC}(f_{\#}) = {\sf Dom}(f) \, \# $.

\begin{pro}\hspace{-.06in}{\bf .}  

\noindent (1) \ For any $L \subseteq \{0,1\}^*$, $L \#$ is a prefix code 
in $\{0,1,\#\}^*$. 

\smallskip

\noindent (2) \ $L$ is in {\sf P} \ iff \ $L \#$ is in {\sf P}.

\smallskip

\noindent (3) \ For any partial function $f: \{0,1\}^* \to \{0,1\}^*$,
 \, $f_{\#}$ is a right ideal morphism of $\{0,1,\#\}^*$.

\smallskip

\noindent 
(4) \ $f \in {\sf fP}$ \ iff \ \, $f_{\#} \in {\cal RM}^{\sf P}_3$.
\end{pro}
{\bf Proof.} This is straightforward.  \ \ \ $\Box$

\medskip

\bigskip

%%%%%%%%%%%%%%%%%%%
\noindent {\bf Coding from three letters to two letters:}
  
\medskip

\noindent We consider the following encoding from the 3-letter alphabet 
$\{0, 1, \#\}$ to the 2-letter alphabet $\{0,1\}$. 

\smallskip

 \ \ \ \ \ \ ${\sf code}: \, \{0,1,\#\} \, \to \, \{00, 01, 11\}$ 
 \ \ is defined by

\smallskip

 \ \ \ \ \ \ ${\sf code}(0) = 00$, \ \ \ ${\sf code}(1) = 01$, 
 \ \ \ ${\sf code}(\#) = 11$.

\smallskip

\noindent For a word $w \in \{0,1,\#\}^*$, ${\sf code}(w)$ is the 
concatenation of the encodings of the letters in $w$. 

The choice of this code is somewhat arbitrary; e.g., we could 
have picked the encoding 
$c$ from $\{0,1,\#\}$ onto the maximal prefix code $\{00, 01, 1\}$, defined 
by $c(0) = 00$, $c(1) = 01$, $c(\#) = 1$. 

\begin{defn} \label{f_encoding} 
 \ We define $f^C: \{0,1\}^* \to \{0,1\}^*$ by letting 
 \, ${\sf Dom}(f^C) \ = \ {\sf code}({\sf Dom}(f) \, \#) \ \{0, 1\}^*$, and 

\smallskip

 \ \ \ \ \ \   
$f^C({\sf code}(x \#) \, v) \ = \ {\sf code}(f(x) \, \#) \ v$,

\smallskip

\noindent for all $x \in {\sf Dom}(f)$ $( \subseteq \{0,1\}^*)$,
and all $v \in \{0,1\}^*$.
We call $f^C$ the {\em encoding of $f$}. 
\end{defn}

\begin{pro}.

\noindent (1) \ For any $L \subseteq \{0,1\}^*$, \, ${\sf code}(L \#)$ is a 
prefix code.

\smallskip

\noindent (2) \ $L$ is in {\sf P}
 \ iff \ ${\sf code}(L \#)$ is in {\sf P}.

\smallskip

\noindent (3) \ For any partial function $f: \{0,1\}^* \to \{0,1\}^*$, 
 \, $f^C$ is a right ideal morphism of $\{0,1\}^*$.

\smallskip

\noindent (4) \ $f \in {\sf fP}$ 
 \ iff \ $f^C \in {\cal RM}_2^{\sf P}$.
\end{pro}
{\bf Proof.} This is straightforward. \ \ \ $\Box$

\begin{pro}.

\noindent (1) \ The transformations 
$f \in {\sf fP} \longmapsto f_{\#} \in {\cal RM}_3^{\sf P}$ and  
$f \in {\sf fP} \longmapsto f^C \in {\cal RM}_2^{\sf P}$ are injective
total homomorphisms from {\sf fP} into ${\cal RM}_3^{\sf P}$, respectively
${\cal RM}_2^{\sf P}$.

\smallskip

\noindent (2) \ $f$ is regular in {\sf fP} \ iff \ $f_{\#}$ is regular in
${\cal RM}^{\sf P}_3$ \ iff \ $f^C$ is regular in ${\cal RM}_2^{\sf P}$.

\smallskip

\noindent (3) \ There are one-to-one correspondences between the inverses
of $f$ in {\sf fP}, the inverses of $f_{\#}$ in ${\cal RM}^{\sf P}_3$, and
the inverses of $f^C$ in ${\cal RM}_2^{\sf P}$.
\end{pro} 
{\bf Proof.} \ (1) is straightforward, and (2) follows from injectiveness 
and from the fact that the homomorphic image of an inverse is an inverse. 

(3) Let $G \in {\cal RM}^{\sf P}_3$ be such that 
$f_{\#} \circ G \circ f_{\#} = f_{\#}$; i.e., 
$f_{\#}(G(f(x) \# w)) = f(x) \# w$, for all $x \in \{0,1\}^*$ and 
$w \in \{0,1,\#\}^*$. Since $f_{\#}(G(f(x) \# w))$ ($= f(x) \# w$) contains 
$\#$, $G(f(x) \# w)$ is of the form $G(f(x) \# w) = z \# v$, for some
$z \in \{0,1\}^*$ and $v \in \{0,1,\#\}^*$. Hence $f_{\#}(G(f(x) \# w)) = $
$f_{\#}(z \# v) = f(z) \# v =  f(x) \# w$, so $f(z) = f(x)$ and $v = w$. So,
$G(f(x) \# w) = z \# w$ for some $z \in f^{-1}f(x)$.
Thus there exists a function $g: \{0,1\}^* \to \{0,1\}^*$ such that
$G(y\#) = g(y) \#$ for all $y \in {\sf Im}(f)$; then $G(y\# w) = g(y) \# w$,
for all $w \in \{0,1,\#\}^*$. Hence $g$ is an inverse 
of $f$.  Moreover, $g$ is clearly in {\sf fP} if $G$ is in 
${\cal RM}^{\sf P}_3$. 

Let $H \in {\cal RM}^{\sf P}_2$ be such that $f^C \circ H \circ f^C = f^C$;
i.e., $f^C(H({\sf code}(f(x) \#) \, v)) = {\sf code}(f(x) \#) \, v$, for all 
$x, v \in \{0,1\}^*$. Since $f^C$ outputs ${\sf code}(f(x) \#) \, v$ 
on input $H({\sf code}(f(x) \#) \, v)$, the definition of $f^C$ implies 
that for all $v \in \{0,1\}^*:$ \,  $H({\sf code}(f(x) \#) \ v)$ is of the 
form ${\sf code}(z \#) \, v$ for some $z \in \{0,1\}^*$ with $f(z) = f(x)$. 
Hence there exists a function $h: \{0,1\}^* \to \{0,1\}^*$ such that 
$H({\sf code}(y \#)) =  {\sf code}(h(y) \#)$ for all $y$; then  
$H({\sf code}(y \#) \, v) =  {\sf code}(h(y) \#) \, v$, for all $y,v$.
Hence $h$ is an inverse of $f$.  
Moreover, $h$ is clearly in {\sf fP} if $H$ is in ${\cal RM}_2^{\sf P}$. 
 \ \ \ $\Box$

\bigskip

\noindent We will show in Section 5 that the encoding $f \mapsto f^C$ 
corresponds to an ``inversive reduction''.

\bigskip

\noindent In summary we have the following 
{\bf relation between {\sf fP} and ${\cal RM}_2^{\sf P}$:}

$${\sf fP} \ \stackrel{ \, C}{\hookrightarrow} \ {\cal RM}_2^{\sf P} 
 \ \hookrightarrow \ [ {\sf id} ]_{_{{\cal J}({\sf fP})}}^0 \ \hookrightarrow 
 \ {\sf fP} .$$

\noindent Here $[ {\sf id} ]^0_{_{{\cal J}({\sf fP})}}$ is the 
${\cal J}$-class of the identity {\sf id} of {\sf fP}, together with the 
zero element (i.e., it is the {\it Rees quotient} of the ${\cal J}$-class
of {\sf id} in {\sf fP}).
The embedding into $[ {\sf id} ]^0_{_{{\cal J}({\sf fP})}}$ holds because 
${\cal RM}_2^{\sf P}$ is ${\cal J}^0$-simple.

\bigskip

\section{\bf Evaluation maps}

\bigskip

\noindent The {\em Turing machine evaluation function}
${\sf eval}_{_{\sf TM}}$ is the input-output function of a universal Turing
machine; it has the form
 \ ${\sf eval}_{_{\sf TM}}(w, x) \ = \ \phi_w(x)$,
where $\phi_w$ is the input-output function described by the word 
(``program'') $w$. 
(Recall that by ``function'' we always mean partial function.)
Similarly, there is an {\em evaluation function for acyclic circuits}, 
${\sf eval}_{\sf circ}(C, x) = f_C(x)$, where $f_C$ is the input-output map 
of the circuit $C$. Here we will only consider length-preserving circuits, 
i.e., $|f_C(x)| = |x|$. 
We also identify the circuit with a bitstring that describes the circuit.
The map ${\sf eval}_{\sf circ}$ is polynomial-time computable, but not 
polynomially balanced (since the size of input component $C$ is not bounded 
in terms of the output length $|f_C(x)|$).

Levin \cite{Levin1W} noted that functions of the form 

\smallskip

\hspace{1in}  ${\sf ev}(w, x) \ = \ (w, \phi_w(x))$, 

\smallskip

\noindent (under some additional assumptions) are polynomially balanced and 
polynomial-time computable; and he observed that {\sf ev} is a {\it critical 
one-way function} in the following sense: 

\begin{defn}
A function $e \in {\sf fP}$ is {\em critical} (or {\sf fP}{\em -critical}) 
iff the following holds:
One-way functions exist iff the function $e$ is a one-way function.
Similarly, a set $L \in {\sf NP}$ is {\em critical} (or 
{\sf P}{\em -critical})  iff  the following holds: ${\sf P} \neq {\sf NP}$ 
iff $L \not\in {\sf P}$.
\end{defn}
The literature calls these functions ``universal'' one-way functions; 
however, not all critical one-way functions are universal (in the sense of 
universal Turing machines, or other universal computing devices).
Levin's idea of a ``universal'' (critical) one-way function has also been 
used in probabilistic settings for one-way functions (see e.g., 
\cite{Goldreich}). 

To make {\sf ev} polynomial-time computable, additional features have to be
introduced. One approach is to simply build a counter into the program of 
{\sf ev} that stops the computation of ${\sf ev}(w, x)$ after a polynomial
number of steps (for a fixed polynomial). For example, the computation could 
be stopped after a quadratic number of steps, i.e., $c \, (|w| + |x|)^2 + c$ 
steps (for a fixed constant $c \geq 1$); we call this function 
${\sf ev}_{(2)}$. 
There exist other approaches; see e.g., section 2.4 of \cite{Goldreich}, 
or p.\ 178 of \cite{AroraBarak}, where it is proved that ${\sf ev}_{(2)}$ 
is {\sf fP}-critical.

Here is another simple example of a critical one-way function:

\smallskip

\hspace{1in} ${\sf ev}_{\sf circ}(C, x) \ = \ (C, f_C(x))$,

\smallskip

\noindent where $C$ ranges over finite acyclic circuits (more precisely, 
strings that describe finite acyclic circuits), and $f_C$ is the 
input-output map of a circuit $C$. 
The function ${\sf ev}_{\sf circ}$ is in {\sf fP}; it is balanced since 
$|x| \leq |C|$ and $C$ is part of the output. Here we only consider
length-preserving circuits, i.e., $|f_C(x)| = |x|$.
We will prove later that ${\sf ev}_{\sf circ}$ is not only critical, but 
also {\em complete} with respect to a reduction that is appropriate for 
one-way functions. 

A similar example of a critical function is 
$(B, \tau) \longmapsto (B, B(\tau))$, where $B$ ranges over all boolean
formulas (or over all boolean formulas in 3CNF), $\tau$ is any truth-value 
assignment for $B$ (i.e., a bitstring whose length is the number of boolean
variables in $B$), and $B(\tau)$ is the truth-value of $B$ for the 
truth-value assignment $\tau$. More generally we have:

\begin{pro} \label{f_Mcrit} \,   
For a nondeterministic polynomial-time Turing machine $M$, let $f_M$ be 
defined by  

\smallskip

 \ \ \ $f_M(x,s) = x$ \ iff \ $M$, with {\it choice sequence} $s$,
accepts $x$ (and undefined otherwise).

\smallskip

\noindent Then $f_M$ is {\sf fP}-critical iff the language (in {\sf NP})
accepted by $M$ is {\sf P}-critical. 
\end{pro}
{\bf Proof.} We studied the functions $f_M$ in Prop.\ \ref{ImfP} (where we 
used the notation $f_L$). We saw that $f_M$ is one-way iff 
${\sf Im}(f_M) \not\in {\sf P}$. Moreover, ${\sf Im}(f_M)$ is the language 
accepted by $M$, and ${\sf Im}(f_M) \in {\sf NP}$. So, $f_M$ is one-way iff 
${\sf P} \neq {\sf NP}$.  
 \ \ \ $\Box$ 

\bigskip

\bigskip

\bigskip

\noindent {\bf Machine model for {\sf fP}}:

\medskip

\noindent Every function in {\sf fP} can be computed by a Turing machine
with a built-in polynomial-time counter, that is used for enforcing 
time-complexity and input balance. 
As usual, to say that time-complexity or balance functions are 
``polynomial'' means that they have polynomial upper-bounds.
More precisely, we will describe every polynomial-time multi-tape Turing 
machine $M$ by a program $v$ (which consists of the list of transitions of 
the Turing machine, as well as its start and accept states), and a 
polynomial $p$ such that $p(n)$ is an upper-bound on the time-complexity 
and the input balance of $M$ on all inputs of length $\le n$. Since we only 
require polynomial {\it upper-bounds}, we can take $p$ of the form 
$p(n) = a \, n^k + a$, where $k, a$ are positive integers.
So $p$ is determined by two integers (stored as bitstrings). 
We do not need to assume anything about the time-complexity of the Turing 
machine with program $v$ (and in general it is undecidable whether $v$
has polynomial-time complexity); instead, we want to consider pairs $(v,p)$ 
where $v$ is a Turing machine program, and $p$ is a polynomial (given by 
two integers $k,a$). 
Based on pairs $(v,p)$ we define the following: A partial function $f$ 
is computed by $(v,p)$ iff for all $x \in A^*$, $f(x)$ is computed by the  
Turing machine with program $v$ in time $\le a \, |x|^k + a = p(|x|)$ and 
input balance $\le p(|f(x)|)$; when $f(x)$ is undefined then the program
either gives no output, or violates the time bound or the input balance
bound.
In this way {\sf fP} can be recursively enumerated by pairs $(v,p)$. 
For {\sf P} and {\sf NP} this (or a similar idea) goes back to the work of 
Hartmanis, Lewis, Stearns, and others in the 1970's; compare with the 
generic {\sf NP}-complete problem in \cite{Hartmanis}, the proof of the 
complexity hierarchy theorems in chapter 12 in \cite{HU}, and the section 
on critical (``universal'') one-way functions in \cite{Goldreich}.  

However, pairs $(v,p)$ do not form a machine model, being hybrids
consisting of a machine and two numbers.
In order to obtain a machine model for {\sf fP} we take a Turing machine 
with program $v$, and add an extra tape that will be used as a counter. 
We assume that every tape has a left endmarker.  On input $x$, a Turing 
machine with counter first computes $p(|x|) = a \, |x|^k + a$, and 
moves the head of the counter tape $p(|x|)$ positions to the right. 
After the counter has been prepared (and the head on the input $x$ has been 
moved back to the left end), the Turing machine executes program $v$ on the 
other tapes, while in each transition the head on the counter tape moves 
left by one position. If the counter head gets back to the left endmarker
(``it triggers the counter''), the Turing machine stops and rejects (and 
produces no output). If the machine halts before triggering the counter, 
the counter has no effect on the result of program $v$ on input $x$. After 
this, if there is an output $y$ the Turing machine with counter checks the 
input balance: If $|y| \ge |x|$ the balance condition obviously holds, so
$y$ is the final output. If $|y| < |x|$, the machine computes $p(|y|)$ 
($< p(|x|)$); if $|x| > p(|y|)$ the machine rejects (and produces no final 
output); otherwise, $y$ is the final output.    

\smallskip

In order to mark off space of length $p(|x|)$ on the counter tape, we need 
an algorithm for computing $p(|x|) = a \, |x|^k + a$, and we will look at 
the time-complexity of this algorithm.  Recall that the bitlength of a 
positive integer $n$ is $\lfloor \log_2 n \rfloor + 1$ (in unsigned binary
representation). 

\smallskip

\noindent 
(1) First, we compute $|x|$ in binary, by repeatedly dividing $|x|$ by 2, 
using two tapes: On one tape, we start with a length $n = |x|$, 
then a mod-2 counter produces $\lfloor n/2 \rfloor$ on the 2nd tape and 
records the remainder (0 or 1) on a 3rd tape; then a mod-2 counter computes 
half of the 2nd tape and writes it on the 1st tape, and records the 
remainder on tape 3, etc. This takes time 
 \ $\le \sum_{i=0}^{\lfloor \log_2 |x| \rfloor} |x|/2^i$ 
$ \ = \ 2^{\lfloor \log_2 |x| \rfloor + 1} - 1 \ < \ 2 \, |x|$.

\smallskip

\noindent 
(2) We compute $a \, |x|^k + a$ in binary, using $k$ multiplications 
and one addition. This takes time 

\smallskip

$\le \ c \, (\log_2 |x|)^2 + 2 c (\log_2 |x|)^2 + \ \ldots \ + $
$ (k-1) c (\log_2 |x|)^2$ $+$ $c k \log_2 |x| \ \log_2 a$

 \ \ \ \ \ $+$ $ \ c\, \max\{k \log_2 |x|, \log_2 a\}$ 

\smallskip

$< \ c \, k^2 \, (\log_2 |x|)^2 + c \, k \ \log_2 a + c \, k \ \log_2 |x|$ 
 \ \ \ \ \ (when \ $k \log_2 |x| \ge \log_2 a$, i.e., $|x| \ge a^{1/k}$)

\smallskip

$< \ c \ k^2 \ \log_2 a \ (\log_2 |x|)^2$,

\smallskip

\noindent where $c$ is a positive constant that depends on the details of 
the multiplication and addition algorithms. Two integers $n_1, n_2$ (in 
binary) can be multiplied in time $\le c \, \log_2 n_1 \, \log_2 n_2$, and 
added in time $\le c \, \max\{\log_2 n_1, \ \log_2 n_2\}$. The bitlength 
of the product $n_1 n_2$ is $\le$ $\lfloor \log_2 n_1 \rfloor$ $+$ 
$\lfloor \log_2 n_2\rfloor + 2$. 
For the last expression in the calculation above, we have 
 \ $c \ k^2 \cdot \log_2 a \cdot (\log_2 |x|)^2 \ < \ |x|$ 
 \ for large enough $|x|$, i.e., when $c_p \le |x|$ \, (where $c_p$ is 
a positive integer depending on  $p$).  

Remark (concrete upper-bound on $c_p$): We define $c_p$ to be the smallest 
number $N$ such that for all $x \in A^{\ge N}$, 
the time to prepare the counter is $\le |x|$.  We saw that this time is
$\le |x|$ when \ $|x| \ge a^{1/k}$ \ and
 \ $|x| \ge c\, k^2 \, \log_2 a \ (\log_2 |x|)^2$.
One proves easily that $n \ge (\log_2 n)^2$ for all $n \ge 16$.
So, $|x|  = |x|^{1/2} \cdot |x|^{1/2}$
$ \ge c\, k^2 \, \log_2 a \cdot (\log_2 |x|)^2$ \ is implied by
 \ $|x|^{1/2} \ge c\, k^2 \, \log_2 a$ \ and \ $|x|^{1/2} \ge 16$. 
Thus we have: \ $c_p \ \le \ $
$\max\{256, \ \ c^2 \, k^4 \, (\log_2 a)^2, \ \ a^{1/k} \}$.

\smallskip

\noindent 
(3) We mark off space of length $p(|x|)$ by converting the binary 
representation of $p(|x|)$ to ``unary'': This is done by the Horner scheme 
with repeated doubling (where the doubling is done by using two tapes, 
and writing two spaces on the 2nd tape for each space on the 1st tape). 
This takes time 
 \ $\sum_{i=0}^{\lfloor \log_2 p(|x|) \rfloor + 1} 2^i \ < \ 4 \, p(|x|)$. 

\smallskip

Thus, the total time used to prepare the counter on input $x$ is 
$< \ 2 \, |x| + |x| + 4 \, p(|x|) \, \le \, 7 \, p(|x|)$, when 
$|x| \ge c_p$. 

For inputs $x$ with $|x| < c_p$, the counter will also receive space 
$p(|x|)$, but the time used for this could be more than $7 \, p(|x|)$.
We can remove the exception of the finitely many inputs of length $< c_p$ 
as follows. For these inputs we let the Turing machine operate as a
finite-state machine (without using the work tapes or any counter); for 
such an input $x$, the time to set up the counter will then be 
 \ $\le \, |x| + p(|x|)$ ($< 7 \, p(|x|)$). 

After execution of program $v$ on input $x$ for time $\le p(|x|)$, if
there is an output $y$ so far, the input balance is checked. 
If $|y| \ge |x|$ input balance holds automatically; checking whether
$|y| \ge |x|$ takes time $\le |x|$. If $|x| > |y|$,
we compute $p(|y|)$ ($< p(|x|)$) in binary, in the same ways as in steps 
(1) and (2) of the counter-tape preparation above. 
This takes time \, $\le \, 2 \, |y| + |y| \, < \, 3 \, |x|$, when 
$|y| \ge c_p$. The time needed to compare the binary representations of 
$|x|$ and $|y|$ (of length $\le \lfloor \log_2 |x| \rfloor + 1$) is 
absorbed in the time to calculate $p(|y|)$. 
So, the time for checking the input balance is $< |x| + 3 \, |x|$
$<$ $4 \, p(|x|)$.

\smallskip

So far we have obtained a machine, described by $(v,p)$, whose 
time-complexity is $\le 7 \, p(.) + p(.) + 4\, p(.) = 12 \ p(.)$ and whose 
input balance is $\le p(.)$.
Because of the preparation of the counter, the time-complexity of the new 
machine is always larger than $p(.)$; therefore we will further modify the 
construction, as follows. Let $p(n) = a \, (n^k + 1)$; we will assume from 
now on that $a \ge 12$. Let $p'(n) = (a - a\%12) \cdot (n^k + 1)$, where 
$a\%12$ is the remainder of the division of $a$ by 12. 
So, $a - a\%12 \ge 12$ and $a - a\%12$ is a multiple of 12.   
 \ (1) Instead of marking off space of length $p(|x|)$, the new machine marks 
off length $\frac{1}{12} \, p'(|x|)$, in time $\leq \frac{7}{12} \, p'(|x|)$.
 \ (2) It executes the program $v$ for time $\le \frac{1}{12} \, p'(|x|)$, 
using the marked-off counter.
(3) It checks the input balance in time $\le \frac{4}{12} \, p'(|x|)$. 
 The total time of the modified machine is then 
$\le p'(|x|) \le p(|x|)$, and the input balance is $\le p'(|x|) \le p(|x|)$.

So for every Turing machine program $v$ and any polynomial $p$ this 
modified machine computes a function in {\sf fP}. Conversely, if $v$ is a 
program with polynomial time for a function $f \in {\sf fP}$, then the 
modified program with pair $(v,p)$ correctly computes $f$ if 
$\frac{1}{12} \, p'(.)$ is larger than the time (and balance) that $v$ uses 
on all inputs; if $f$ belongs to {\sf fP} then a polynomial $p$ such that
$\frac{1}{12} \, p'(.)$ is large enough for bounding the time and the input 
balance, exists. 
Such a modified machine will be called {\em Turing machine with polynomial
counter}. A program for such a machine consists of a pair $(v,p)$ and an
extra program for preparing the counter and checking input balance; let's
call that extra program $u_p$ (it only depends on $p$, and not on $v$). 
The triple $(v, p, u_p)$, or more precisely, the word 
${\sf code}(v \# k \# a \# u_p)$ with the numbers $k,a$ written in 
binary, will be called a {\it polynomial program}.  Thus, Turing machines 
with polynomial counter are a machine model for {\sf fP}.
So we have proved:

\begin{pro} \label{Polyn_TM}
There exists a class of modified Turing machines, called {\em Turing machine 
with polynomial counter}, with the following properties:
For every $f \in {\sf fP}$ there exists a Turing machine with polynomial 
counter that computes $f$; and for every Turing machine with polynomial
counter, the input-output function belongs to {\sf fP}.  
  \ \ \ $\Box$
\end{pro} 
{\bf Notation:} \ A polynomial program 
$w = {\sf code}(v \# k \# a \# u_p)$, based on a Turing machine program $v$ 
and a polynomial $p$ (with $p(n) = a \, n^k + a$), will be denoted by 
$\langle v,p \rangle$. The polynomial $p$ that appears in $w$ will often be 
denoted by $p_w$. The function computed by a polynomial program 
$w = \langle v,p \rangle$ will be denoted by $\phi_w$ ($\in {\sf fP}$).

\bigskip

\noindent {\bf Evaluation maps for {\sf fP}}:

\medskip

At first we consider a function ${\sf ev}_{\sf poly}$ defined by 
${\sf ev}_{\sf poly}(w, x) = (w, \phi_w(x))$, 
where $w = \langle v,p \rangle$ is any polynomial program.
But ${\sf ev}_{\sf poly}$ is {\em not} in {\sf fP}. Indeed, the output 
length (and hence the time-complexity) of  ${\sf ev}_{\sf poly}$ on input 
$(w,x)$ is equal to $p(|x|)$ (in infinitely many cases, when $p$ is a tight 
upper-bound); as $w$ varies, the degree of $p$ is unboundedly large, hence 
the time-complexity of ${\sf ev}_{\sf poly}$ has no polynomial upper-bound. 
We will nevertheless be able to build {\sf fP}-critical functions. For a 
fixed polynomial $q$ (of the form $q(n) = a \, n^k + a$), let

\medskip

\hspace{1in}   ${\sf fP}^{q} \ = \ \{ \phi_w \in {\sf fP} : \ p_w \le q\}$.

\medskip

\noindent More explicitly, $p_w \le q$ means that the polynomial program $w$
has time-complexity $\le p_w(|x|) \le q(|x|)$ and input-balance 
$|x| \leq p_w(|\phi_w(x)|) \le q(|\phi_w(x)|)$, for all 
$x \in {\sf Dom}(\phi_w)$.

In general, for polynomials $q_1, q_2$ we say $q_1 \le q_2$ iff for all 
non-negative integers $n$: $q_1(n) \le q_2(n)$.
Interestingly, for polynomials of the form $q_i(n) = a_i \, n^{k_i} + a_i$ 
we have: 
$q_1 \le q_2$ \ iff \ $k_1 \le k_2$ and $a_1 \le a_2$.  
Hence it is easy to check whether $p_w \le q$, given $w$ and the two numbers
that determine $q$.  

\smallskip

A polynomial program $w$ such that $p_w \le q$ (for a fixed polynomial $q$) 
is called a {\em $q$-polynomial program}. We define ${\sf ev}_{q}$ by 

\smallskip

\hspace{1in}  ${\sf ev}_{q}(w, x) = (w, \phi_w(x))$,

\smallskip

\noindent where $w$ is any $q$-polynomial program. 
The function ${\sf ev}_{q}$ above has two input and two output strings. 
To make ${\sf ev}_{q}$ fit into our framework of functions with one input 
and one output string we encode ${\sf ev}_{q}$ as 
${\sf ev}_{q}^C: \{0,1\}^* \to \{0,1\}^*$ where for all 
$w, x \in \{0,1\}^*$ such that $w$ is a $q$-polynomial program, 

\medskip

\hspace{1in}   ${\sf ev}_{q}^C\big({\sf code}(w) \, 11 \, x\big) = $
        ${\sf code}(w) \, 11 \, \phi_w(x)$.

\medskip

\noindent From now on we will call ${\sf ev}_{q}^C$ an {\it ``evaluation 
map''}.  We observe that in the special case where $\phi_w$ (for a fixed 
$w$) is a right ideal morphism, the function 

\medskip

\hspace{1in} ${\sf ev}_{q}^C\big({\sf code}(w) \, 11 \, \cdot\big)$:
$ \ \ x \ \longmapsto \ \ {\sf code}(w) \, 11 \, \phi_w(x)$ 

\medskip

\noindent is also a right ideal morphism.

\bigskip

\noindent {\bf Criticality of ${\sf ev}_{q}^C$:}

\medskip

\noindent For any fixed word $v \in \{0,1\}^*$ we define the prepending 
map

\smallskip

\hspace{1in}  $\pi_v: x \in \{0,1\}^* \longmapsto v \, x$ ;

\smallskip

\noindent and for any fixed positive integer $k$ we define 

\smallskip

\hspace{1in} $\pi_k': z \, x \in \{0,1\}^* \longmapsto x$, \ where $|z| = k$
 \ (with $\pi_k'(t)$ undefined when $|t| < k$). 

\smallskip

\noindent Clearly, $\pi_v, \pi_k' \in {\cal RM}_2^{\sf P}$, and we have 
 \ $\pi'_{|v|} \circ \pi_v = {\sf id}_{A^*}$ 
 \ (i.e., $\pi'_{|v|}$ is a left inverse of $\pi_v$).

We observe that $\pi_v$ can be written as a composite of the
maps $\pi_0$ and $\pi_1$, for any $v \in \{0,1\}^*$. Similarly, $\pi_k'$ 
is the $k$th power of $\pi_1'$.

\begin{pro} \label{ev_q} 
 \ Let $q$ be any polynomial such that for all $n \ge 0$, $q(n) > cn+c$
(where $c > 1$ is a constant).  
Then ${\sf ev}_{q}^C$ belongs to {\sf fP}, and ${\sf ev}_{q}^C$ 
is a one-way function if one-way functions exist. 
\end{pro}
{\bf Proof.} We saw that testing whether $p_w \le q$ is easy for 
polynomials of the form that we consider. 
By reviewing the workings of a universal Turing machine (e.g., 
in \cite{HU}) we see that the time-complexity of ${\sf ev}_{q}(w,x)$ is \, 
$\leq \ c_0 \ |w| \cdot p_w(|x|)^2$ (when $p_w$ is at least linear); here, 
$c_0 \ge 1$ is a constant (independent of $x$ and $w$). 
The factor $p_w(|x|)^2$ comes the fact that Turing machines can have any 
number of tapes, whereas a Turing machine for ${\sf ev}_{q}$ has a fixed 
number of tapes; any number of tapes can be converted to one tape, but the
complexity increases by a square (the more efficient Hennie-Stearns 
construction converts any number of tapes to two tapes, with a complexity
increase from $T$ to $T \log T$). The universal Turing machine 
simulates each transition of program $w$ (modified into a 1-tape Turing 
machine) using $\le c_1 \, |w|$ steps (for a constant $c_1 \ge 1$). 

For input balance: When $|x| \le |\phi_w(x)|$ we also have 
$|w| + |x| \le |w| + |\phi_w(x)|$ so balance is automatic.
When $|x| > |\phi_w(x)|$ then (since $p_w$ bounds the input balance of
$\phi_w$), the input-length satisfies \,  
$|w| + |x| \leq |w| + p_w(|\phi_w(x)|) \ \leq \ |w| + q(|\phi_w(x)|)$
$\le$ $q(|w| + |\phi_w(x)|)$. 
Hence ${\sf ev}_{q}$, and similarly ${\sf ev}_{q}^C$, belongs to {\sf fP}.

Criticality: 
If the function ${\sf ev}_{q}$ has an inverse $e'_q \in {\sf fP}$, then 
 \, ${\sf ev}_{q} \circ e'_q \circ {\sf ev}_{q}(w,x) = (w, \phi_w(x)$.
Hence, for any function $\phi_w \in {\sf fP}^{q}$ with a fixed program $w$ 
we have:  
$\phi_w \circ \pi_{|w|}' \circ e'_q \circ \pi_w \circ \phi_w = \phi_w$,
where $\pi_{|w|}': (w,v) \mapsto v$, and $\pi_w: y \mapsto (w,y)$. 
For a fixed $w$ we have $\pi_{|w|}', \pi_w \in {\sf fP}$, so \, 
$\pi_{|w|}' \circ e'_q \circ \pi_w$ \, is an inverse of $\phi_w$. Hence, if
${\sf ev}_{q}$ is not one-way, no function in ${\sf fP}^{q}$ is one-way. 
But there are functions (e.g., some $f_M$ as seen in Prop.\ \ref{f_Mcrit}) 
that are {\sf fP}-critical, even when $q$ is linear (since by padding 
arguments one can obtain {\sf NP}-complete languages with nondeterministic
linear time-complexity).  So some {\sf fP}-critical $f_M$ would have an 
inverse in {\sf fP}. The same proof is easily adapted to ${\sf ev}_{q}^C$. 
 \ \ \ $\Box$ 

\bigskip

\noindent {\bf Finite generation of {\sf fP}:}

\medskip

\noindent We will show that ${\sf ev}_{q}^C$ can be used to simulate
universal evaluation maps, and to prove that {\sf fP} is finitely generated. 
This is based on the universality of ${\sf ev}_{q_2}^C$ for 
${\sf fP}^{q_2}$, combined with a padding argument; here $q_2$ is a
polynomial of degree 2, with $q_2(n) \ge  c\, n^2 +c$ for a constant 
$c \ge 12$.  (We chose 12 in view of the reasoning before Prop.\ 
\ref{Polyn_TM}.)
We need some auxiliary functions first.

\medskip

\noindent We define the {\em expansion (or padding) map}, first as a 
multi-variable function for simplicity:

\smallskip

\hspace{.6in} ${\sf expand}(w, x) \ = \ $
    $({\sf e}(w), \ (0^{4 \, |x|^2 + 7 \, |x| + 2}, \ x))$

\smallskip

\noindent where ${\sf e}(w)$ is such that \    

\smallskip

\hspace{.6in} 
$\phi_{{\sf e}(w)}(0^k, \ x) \ = \ (0^k, \ \phi_w(x))$, \ for all $k \geq 0$.

\smallskip

\noindent The program ${\sf e}(w)$ is easily obtained from the program 
$w$, since it just processes the padding in front of the input and in front 
of the output of $\phi_w$, and acts as $\phi_w$ on $x$; and the time and
balance polynomial is decreased due to padding.  
As a one-variable function, {\sf expand} is defined by

\medskip

\hspace{.6in}  ${\sf expand}\big({\sf code}(w) \ 11 \ x\big) \ = \ 
   {\sf code}({\sf ex}(w)) \ 11 \ 0^{4 \, |x|^2 + 7 \, |x| + 2} \ 11 \ x$,

\medskip

\noindent where for one-variable functions, the program ${\sf ex}(w)$ 
is such that \    

\medskip

\hspace{.6in} $\phi_{{\sf ex}(w)}(0^k \ 11 \ x) \ = \ 0^k \ 11 \ \phi_w(x)$.

\medskip

\noindent Again, ${\sf ex}(w)$ is a slight modification of 
$w = {\sf code}(v \# k \# a \# u_p)$ (following the notation for the
machine model for {\sf fP}), to allow inputs and outputs with padding, and 
to readjust the complexity and balance polynomial. 
More precisely, ${\sf ex}(w)$ is of the form 
${\sf code}(r \# v \# \lceil k/2 \rceil \# a_e \# u_{p_e})$, where $r$ is 
a preprocessing subprogram by which the prefix $0^k 11$ of the input is 
simply copied to the output; at the end of execution of $r$, the state and 
head-positions are the start state and start positions of the subprogram 
$w$. The appropriate complexity and balance polynomial stored in 
${\sf ex}(w)$ is 

\medskip

 \ \ \ $p_e(n) = a_e \, n^{\lceil k/2 \rceil} + a_e$, \ with \ 
$a_e \ = \ \max\{12, \ \lceil a/2^k \rceil +1\}$.  

\medskip

\noindent Indeed, if $m = |x|$, an input 
$0^{4 \, |x|^2 + 7 \, |x| + 2} \ 11 \ x$ of $\phi_{{\sf ex}(w)}$ has 
length $i = 4 \, m^2 + 8 \, m + 4$. So, $m = \sqrt{i}/2 - 1$.  
Let $a \, (m^k + 1)$ be the polynomial of program $w$.
The complexity of $\phi_{{\sf ex}(w)}$ on its input 
is $4 \, m^2 + 7 \, m + 4$ (for reading the part $0^* 11$ of the input), 
plus \, $a \, (m^k + 1)$ (for using $x$ and computing $\phi_w(x)$).  
So in terms of its input length $i$, the complexity of 
$\phi_{{\sf ex}(w)}$ is \ $< \ i + a \, (m^k + 1)$
$\le \ i^{\lceil k/2 \rceil} \ + \ a \, ((\sqrt{i}/2 - 1)^k + 1)$ 
$\le \ i^{\lceil k/2 \rceil} \ + \ a/2^k \, i^{k/2}$; 
the last step uses the fact that $(z - 1)^m \le z^m - 1$ for all 
$z \ge 0, \ m \ge 1$.  Hence the complexity of $\phi_{{\sf ex}(w)}$ is 
$< \ (a/2^k + 1) \, i^{\lceil k/2 \rceil}$.  
For the input balance of $\phi_{{\sf ex}(w)}$ we have: The input-length
is bounded by twice the output-length.
Indeed, the input length is $i = 4 \, |x|^2 + 8 \, |x| + 4$
$<$ $2 \cdot |0^{4 \, |x|^2 + 7 \, |x| + 2} \ 11|$
$<$ $2 \cdot |0^{4 \, |x|^2 + 7 \, |x| + 2} \, 11 \, \phi_w(x)|$ $=$ 
$2 \cdot |\phi_{{\sf ex}(w)}(0^{4 \, |x|^2 + 7 \, |x| + 2} \, 11 \, x)|$.
Moreover, we want $a_e$ to stay $\ge 12$ (in view of the reasoning before 
Prop.\ \ref{Polyn_TM}). 

\smallskip

In order to achieve an arbitrarily large polynomial amount of padding we 
iterate the quadratic padding operation. Therefore we define a {\em repeated 
expansion (or re-padding)} map, first as a two-variable function:

\medskip

\hspace{.6in} ${\sf reexpand}(u, (0^h, x)) \ = \ $
       $({\sf e}(u), (0^{4 \, h^2 + 8 \, h + 2}, x))$, \ for all $h > 0$, 

\medskip

\noindent where ${\sf e}(.)$ is as above. As a one-variable function, 

\medskip

\hspace{.6in} ${\sf reexpand}\big({\sf code}(u) \ 11 \ 0^h \ 11 \ x \big)$ 
$ \ = \ $ 
${\sf code}({\sf ex}(u)) \ 11 \ 0^{4 \, h^2 + 8 \, h + 2} \ 11 \ x$, 
         \ for any $h \geq 0$,

\medskip

\noindent with ${\sf ex}(.)$ as in {\sf expand} above.

\medskip

We also introduce a {\em contraction (or unpadding) map}, which is a partial
left inverse of {\sf expand}. We define {\sf contr} first as a multi-variable 
function:

\medskip

\hspace{.6in} 
$ {\sf contr}(w, \ (0^h, \ y)) \ = \ ({\sf c}(w), \ y)$, 
  \ \ if $h \leq 4 \, |y|^2 + 7 \, |y| + 2$ \ \ (undefined otherwise).  
 
\medskip

\noindent As a one-variable function, {\sf contr} is defined by

\medskip 

\hspace{.6in} 
${\sf contr}\big({\sf code}(w) \ 11 \ 0^h \ 11 \ y \big) \ = \  $
${\sf code}({\sf co}(w))\ 11 \ y$, \ \ if $h \leq 4\, |y|^2 + 7\, |y| + 2$ 
  
\hspace{.6in} (undefined otherwise). 

\medskip

\noindent The program transformations ${\sf c}(.)$ and ${\sf co}(.)$ are
inverses of ${\sf e}(.)$, respectively ${\sf ex}(.)$. So, ${\sf c}(.)$ and
${\sf co}(.)$ erase the prefix $r$ in ${\sf ex}(u)$, and replace the
polynomial $b \, n^h + b$, encoded in ${\sf ex}(u)$, by
$b_c \, n^{2h} + b_c$, where  \, $b_c = (b-1) \, 2^{2h}$.

\smallskip

To invert repeated padding we introduce a {\em repeated contraction (or 
unpadding)} map, first as a multi-variable function. 
Note that if $h = 4 \, k^2 + 8 \, k + 2$ (which is the amount of padding 
introduced by {\sf reexpand}), then $k = \frac{1}{2} \, \sqrt{h+2} -1$.
Therefore, for any $h \ge 0$ we define

\medskip

\hspace{.6in} ${\sf recontr}(u, \ (0^h, \ y)) \ = \ $ 
$({\sf c}(u), \ (0^{\max\{1, \ \lfloor \sqrt{h+2}/2 \rfloor - 1\}}, \ y))$
   \ \ (undefined on other inputs).
 
\medskip

\noindent As a one-variable function, {\sf recontr} is defined by

\medskip

\hspace{.6in} 
${\sf recontr}\big({\sf code}(u) \ 11 \ 0^h \ 11 \ y \big)$
 $ \ = \ $ 
${\sf code}({\sf co}(u))\ 11$
    $0^{\max\{1, \ \lfloor \sqrt{h+2}/2 \rfloor - 1\}} \ 11\ y$ 

\hspace{.6in} (undefined on other inputs).

\medskip

The maps {\sf expand}, {\sf reexpand}, {\sf contr}, and {\sf recontr} belong 
to {\sf fP}, and they are regular (they have polynomial-time inverses).

\begin{pro} \label{fPfingen} 
 \ \ {\sf fP} is finitely generated.
\end{pro}
{\bf Proof.} \ We will show that the following is a generating set of 
{\sf fP}: 

\medskip
 
\hspace{1in} 
$\{{\sf expand}, \ {\sf reexpand}, \ {\sf contr}, \ {\sf recontr}, $
$ \ \pi_0, \ \pi_1, \ \pi_1', \ {\sf ev}_{q_2}^C \}$, 

\medskip

\noindent where $q_2$ is the polynomial $q_2(n) = c \, n^2 + c$, with
$c \geq 12$ (the number 12 comes from the discussion before Prop.\
\ref{Polyn_TM}).

\smallskip

\noindent {\sf Remark:} \, The functions {\sf expand}, {\sf reexpand}, 
{\sf contr}, and {\sf recontr} all have quadratic time-complexity and
balance functions, so they can be generated by $\pi_0, \pi_1, \pi_1'$, and
${\sf ev}_{q_2}^C$ (provided that the constant $c$ in $q_2$ is chosen large 
enough, so that ${\sf ev}_{q_2}^C$ can execute {\sf expand}, {\sf reexpand}, 
{\sf contr}, and {\sf recontr}).  Thus \, 

\medskip

\hspace{1in} 
$\{\pi_0, \ \pi_1, \ \pi_1', \ {\sf ev}_{q_2}^C \}$ \, 

\medskip

\noindent is a generating set of {\sf fP}.  
We use the larger generating set since it yields simpler formulas.

\medskip

Let $w$ be a program with polynomial counter and let $m$ be an integer 
upper-bound on $\log_2(a + k)$, where $p_w (n) = a \, n^k + a$. 
We also assume that the program $w = \langle v, p \rangle$ is such that 
for any polynomial $P(.) > p_w(.)$, $\langle v, P \rangle$ also computes 
$\phi_w$; indeed, since $\phi_w \in {\sf fP}$, we can choose $v$ and $p$ so 
that the execution of $v$ by the Turing machine with polynomial counter 
(described by $\langle v, p \rangle$) never triggers the counter; in that 
case, making the counter larger does not change the function.  
Then for all $x \in \{0,1\}^*$,

\bigskip

\noindent $(\star)$ 
\hspace{.5in} $\phi_w(x) \ = \ $
$\pi_{_{2 \, |w'| + 2}}' \circ {\sf contr} \circ {\sf recontr}^{2 \, m}$ 
$\circ$ ${\sf ev}_{q_2}^C$ $\circ$ 
${\sf reexpand}^{2m} \circ {\sf expand}$  $\circ$
$\pi_{_{{\sf code}(w) \, 11}}(x)$ ,

\bigskip

\noindent where $w' = {\sf co}^{2m+1} \circ {\sf ex}^{2m+1}(w)$.
Indeed,

\medskip

 $x \ \ \stackrel{\pi_{_{{\sf code}(w) \, 11}}}{\longmapsto} \ \ $
${\sf code}(w) \ 11 \ x \ \ $
$\stackrel{\sf expand}{\longmapsto} \ \ $
${\sf code}({\sf ex}(w)) \ 11 \ 0^{4 \, |x|^2 + 7\, |x| +2} \ 11 \ x$ 

\medskip

$\stackrel{{\sf reexpand}^{2m}}{\longmapsto} \ \ $ 
${\sf code}({\sf ex}^{2m+1}(w)) \ 11 \ 0^{N_{2m+1}} \ 11 \ x$. 

\medskip

\noindent Here $N_1 = |x|^2 + 7\, |x| +2$, and 
$|0^{N_1} \, 11 \, x| = (2 \, (|x| + 1))^2$;     inductively, 
$N_i = 4\, N_{i-1}^2 + 8\,  N_{i-1} + 2$ \ for $1 < i \leq 2m+1$,
and $|0^{N_i} \, 11| = (2 \, (N_{i-1} + 1))^2$.

\noindent Continuing the calculation, 

\medskip

$\stackrel{{\sf ev}_{q_2}^C}{\longmapsto} \ \ $
${\sf code}({\sf ex}^{2m+1}(w)) \ 11 \ 0^{N_{2m+1}} \ 11 \ \phi_w(x) $

\medskip

$\stackrel{{\sf recontr}^{2m}}{\longmapsto} \ \ $
${\sf code}(w') \ 11 \ 0^{\ell} \ 11 \ \phi_w(x) \ \ $

\medskip

$\stackrel{\sf contr}{\longmapsto} \ \ $
${\sf code}(w') \ 11 \, \phi_w(x) \ \ $
$\stackrel{\pi_{_{2 \, |w'| + 2}}'}{\longmapsto} \ \ \phi_w(x)$

\medskip

\noindent where $w' = {\sf co}^{2m+1} \circ {\sf ex}^{2m+1}(w)$.
We use $2m$ in ${\sf recontr}^{2 m}$ because $\phi_w(x)$ could be much 
shorter than $x$ (but by input balance, $|x| \le p_w(|\phi_w(x)|)$).
As a consequence, we also use $2m$ in ${\sf reexpand}^{2 m}$ in order to have 
equal numbers of program transformations ${\sf ex}(.)$ and ${\sf co}(.)$.   
Note that doing more input padding than necessary does not do any harm; also,
if $w$ contains a polynomial $p_w$ larger than needed for computing 
$\phi_w$, this does not cause a problem (by our assumption on $v$).  
By the choice of $2m$, the value of $\ell$ above is less than 
$4 \, |\phi_w(x)|^2 + 7 \, |\phi_w(x)| + 2$, so {\sf contr} can be applied 
correctly.  

The argument of ${\sf ev}_{q_2}^C$ in the above calculation has length
$> N_{2m+1} + 2 + |x|$, which is much larger than the time it 
takes to simulate the machine with program $w$ on input $x$ (that time is 
$< c_0 |w| \, p_w(|x|)^2$). In fact, by the choice of $m$, the polynomial 
encoded in ${\sf ex}^{2m+1}(w)$ is the linear polynomial $12 \, (n + 1)$ 
(which is $< q_2(n)$).  Hence, ${\sf ev}_{q_2}^C$ works correctly on its 
input in this context.
 \ \ \ $\Box$

\medskip

We saw that {\sf fP} does not have an evaluation map in the same sense as
the Turing evaluation map. However, formula $(\star)$ in the proof of  
Prop.\ \ref{fPfingen} shows that the map ${\sf ev}_{q_2}^C$ simulates every
function in {\sf fP}, in the following sense: 
$f_2$ {\em simulates} $f_1$ (denoted by $f_1 \preccurlyeq f_2$) iff there 
exist $\beta, \alpha \in {\sf fP}$ such that 
$f_1 = \beta \circ f_2 \circ \alpha$; this is discussed further at the 
beginning of Section 5.  
Formula $(\star)$ in Prop.\ \ref{fPfingen} implies:

\begin{pro} \label{evalfP} 
 \ Every function $f \in {\sf fP}$ is simulated by ${\sf ev}^C_{q_2}$.
 \ \ \  \ \ \  $\Box$
\end{pro}  
It follows from this and the definition of simulation that
${\sf ev}^C_{q_2}$ belongs to the $\cal J$-class of ${\sf id}_{A^*}$ in 
{\sf fP}. 

\medskip

Since {\sf fP} is finitely generated we now have two ways of representing 
each element $g \in {\sf fP}$ by a word: 
(1) We have $g = \phi_w$ for some polynomial program $w \in A^*$ (as seen 
in Prop.\ \ref{ev_q}), and 
(2) $g$ can be represented by a string of generators (considering the finite 
set of generators of {\sf fP} as an alphabet). The next proposition 
describes the translation between these two representations.

\begin{pro} \label{compilers}
There are total computable maps $\alpha, \beta$ such that for any word $s$ 
over a finite generating set of {\sf fP}, $\alpha(s)$ is a polynomial 
program for the function given by $s$; and for any polynomial program $u$, 
$\beta(u)$ is a word for $\phi_u$ over the generators of {\sf fP}.

More precisely, let $\Gamma$ be a finite generating set of {\sf fP}. 
For any $s \in \Gamma^*$, let $\Pi s \in {\sf fP}$ be the element of 
{\sf fP} obtained by composing the generators in the sequence $s$. 
There exist total recursive {\em ``compiler maps''} 
$\alpha: \Gamma^* \to \{0,1\}^*$ and $\beta: \{0,1\}^* \to \Gamma^*$ such 
that for all $s \in \Gamma^*$ and all $w \in  \{0,1\}^*$: \, 
$f_{\alpha(s)} = \Pi s$, and $\Pi \beta(w) = \phi_w$.  
\end{pro}

\smallskip

\noindent
{\bf Proof.} The map $\beta$ is given by formula $(\star)$ in the proof 
of Prop.\ \ref{fPfingen}, where a representation over the generators is 
explicitly constructed. When $u$ is not a well-formed polynomial program
we let $\beta(u)$ be a sequence of generators for the empty function.  

Conversely, by composing a sequence of generators, a function in {\sf fP} 
is obtained (note that every sequence $s$ of generators has a finite 
length).  More precisely, if $f_1, f_2 \in {\sf fP}$ have as complexity 
and balance bounds the polynomials 
$q_i(n) = a_i \, n^{k_i} + a_i$ ($i = 1, 2$), then 
$f_2 \circ f_1$ has input balance $\le q_2 \circ q_1(n)$ (obviously), and 
time-complexity $\le q_1(n) + q_2 \circ q_1(n)$. 
Indeed, a polynomial-time program for $f_2 \circ f_1$ is obtained by first 
taking the program for $f_1$ on input $x$, and then applying the program 
for $f_2$ to $f_1(x)$ (in time $\le q_2(|f_1(x)|)$). 
The corresponding polynomial upper-bound is
 \, $q_1(n) + q_2 \circ q_1(n)$ $=$ 
$a_1 \, (n^{k_1} + 1) + a_2 \, a_1^{k_2} \, (n^{k_1} + 1)^{k_2} + a_2$ $<$
$(a_1 + a_2 \, a_1^{k_2})\, (n^{k_1} + 1)^{k_2} + a_2$.
In order to obtain a polynomial upper-bound of the form $a \, n^k + a$, we
use the inequality  

\smallskip

$(n + 1)^j \ \le \ 2^{j-1} \, (n^j + 1)$, \ for all $n \ge 0, \ j \ge 1$. 

\smallskip

\noindent (To prove this inequality apply calculus to the function 
$f(x) = 2^{j-1} (x^j + 1) - (x + 1)^j$.) 
Thus for $f_2 \circ f_1$ we get a complexity and balance upper-bound 

\smallskip

$q(n) = a \, n^{k_1 k_2} + a$, where
 \ $a = a_2 + a_1 + a_2 \, a_1^{k_2} \, 2^{k_2}$.
\smallskip

\noindent This yields an algorithm for obtaining a polynomial program for 
$f_2 \circ f_1$ from polynomial programs for $f_1$ and $f_2$. For a sequence
of generators $s$, this algorithm can be repeated $|s| - 1$ times to yield
a polynomial programs for the sequence $s$ of generators.  
 \ \ \ $\Box$ 

\medskip

\noindent A finite generating set $\Gamma$ for {\sf fP} can be used to 
construct a {\em generator-based evaluation map} for {\sf fP}, defined by
 \ $(s, x) \in \Gamma^* \times A^* \ \longmapsto \ {\sf ev}_{\Gamma}(s,x)$
$=$ $(s, \, (\Pi s)(x))$. 
 However, ${\sf ev}_{\Gamma}$ does not belong to {\sf fP}, for the 
same reasons as we saw at the beginning of Sect.\ 4 for ${\sf ev}_{\sf poly}$.
(But just as for ${\sf ev}_{\sf poly}$ we could restrict ${\sf ev}_{\Gamma}$ 
to a function that belongs to {\sf fP} and that simulates every element 
of {\sf fP}.)

\begin{pro} \label{notFinPres}
 \ {\sf fP} is {\em not} finitely presented. Its word problem is 
co-r.e.,  but not r.e.
\end{pro}
{\bf Proof.} The word problem is co-r.e.: Let $U, V \in \Gamma^*$; using 
Prop.\ \ref{compilers} we effectively find programs $u, v \in A^*$ from 
$U, V$ such that $\phi_u = \Pi U$, $\phi_v = \Pi V$.
If $\phi_u \neq \phi_v$ then by exhaustive search we will find $x$ such that
$\phi_u(x) \neq \phi_v(x)$, thus showing that $U \neq V$ in {\sf fP}.
When $U = V$ in {\sf fP} then this procedure rejects by not halting.

The word problem of {\sf fP} is undecidable, since the equality problem for
languages in {\sf P} can be reduced to this (reducing $L$ to ${\sf id}_L$ 
or to ${\sf id}_{{\sf code}(L\#)}$).
And the equality problem for languages in {\sf P} is undecidable, since the
universality problem of context-free languages can be reduced to the equality
problem for languages in {\sf P}; all context-free languages are in {\sf P}.
The universality problem for context-free language is the question whether
for a given context-free grammar $G$ with terminal alphabet $A$ (with
$|A| \geq 2$), the language generated by $G$ is $A^*$; this problem is
undecidable (see \cite{HU} Thm.\ 8.11).

Since the word problem is co-r.e.\ but undecidable, it is not r.e.  Hence
these finitely generated monoids are not finitely presented (since the word
problem of a finitely presented monoid is r.e.).
 \ \ \ $\Box$

\begin{pro} \ \ {\sf fP} is finitely generated by {\em regular} elements.
\end{pro}
{\bf Proof.} All the listed generators of {\sf fP} are regular, except 
possibly ${\sf ev}_{q_2}^C$.
Let us define a partial function $E_{q_2} \in {\sf fP}^{q_2}$ by \,
 $E_{q_2}(w, x) = (w, \phi_w(x), x)$, when $\phi_w \in {\sf fP}^{q_2}$. 
Obviously, $E_{q_2}$ is not one-way.
But ${\sf ev}_{q_2}$ (as a two-variable function) can be expressed as a
composition of $E_{q_2}$ and the other (regular) generators.
In more detail, ${\sf ev}_{q_2} = \pi'_{q_2} \circ E_{q_2}$, where 
$\pi'_{q_2}(w,z,x) = (w,z)$ if $|z| \le q_2(|x|)$ and $|x| \le q_2(|z|)$. 
So ${\sf ev}_{q_2}^C$ can be replaced by $E_{q_2}^C$ as a generator.
 \ \ \ $\Box$

\begin{pro} \ There are elements of {\sf fP} and of ${\cal RM}_2^{\sf P}$
that are critical (i.e., non-regular if ${\sf P} \neq {\sf NP}$), whose 
product is a non-zero idempotent.
\end{pro}
{\bf Proof.} For $i = 0, 1$, let $e_i \in {\sf fP}$ be defined, as a
two-variable function, by

\smallskip

 \ \ \ \ \ \
$e_i(w,x) \ = \
\left\{
\begin{array}{ll}
(w, \phi_w(x)) & \mbox{if $x \in i\, \{0,1\}^*$,\ 
     $\phi_w(x) \in i\, \{0,1\}^*$,
     and $|\phi_w(x)| = |x|$;} \\
(w, 0^{|x|}) & \mbox{otherwise.}
\end{array}  \right. $

\smallskip

\noindent Then $(e_1 \circ e_0)(w,x) = (w, 0^{|x|})$ for all $(w, x)$,
so $e_1 \circ e_0$ is an idempotent.

To prove that $e_i$ is critical we reduce the satisfiability problem to
the inversion problem of $e_i$.
The reduction for $e_i$ maps a boolean formula $B$ with $n$ variables to
$(b, i^n 1)$, where $b$ is a program such that
$f_b(i \, \tau) = i^n B(\tau)$; i.e., for a truth-value assignment
$\tau \in \{0,1\}^n$, $f_b$ evaluates $B$ on $\tau$, and outputs the 
resulting truth-value, prefixed with $n$ copies of $i$.
If $e_i$ were regular then ${\sf Im}(e_i)$ would be in {\sf P}, by
Prop.\ \ref{ImfP}.
Then satisfiability of $B$ could be checked by a {\sf P}-algorithm,
since  $B$ is satisfiable iff $(b, i^n 1) \in {\sf Im}(e_i)$.
To obtain one-variable functions we can take $e_i^C$.

To prove the proposition for ${\cal RM}_2^{\sf P}$ we define $e_i \in$
${\cal RM}_2^{\sf P}$ for $i = 0, 1$ as follows, first as two-variable
functions:

\smallskip

 \ \ \ \ \ \
$e_i(w,x) \ = \
\left\{
\begin{array}{ll}
(w, \phi_w(x)) & \mbox{if $x \in 0i \, \{0,1\}^*$,
 \ $\phi_w(x) \in 0i \, \{0,1\}^*$, and $|\phi_w(x)| = |x|$;} \\
(w, x)      & \mbox{if $x \in 1 \{0,1\}^*$;} \\
\mbox{undefined} & \mbox{otherwise.}
\end{array}  \right. $

\smallskip

\noindent Then $(e_1 \circ e_0)(w,x) = (w,x)$ when $x \in 1 \, \{0,1\}^*$, 
and $(e_1 \circ e_0)(w,x)$ is undefined otherwise; so $e_1 \circ e_0$ is a
partial identity.
The reduction of the satisfiability problem to the inversion problem of $e_i$
is similar to the case of {\sf fP}.
 \ \ \ $\Box$

\section{Reductions and completeness}

\noindent The usual reduction between partial functions 
$f_1, f_2: A^* \to A^*$ is as follows. 

\begin{defn} \label{simul} \   
 \ $f_1$ {\em is simulated by} $f_2$ (denoted by
$f_1 \preccurlyeq f_2$) \ iff \ there exist polynomial-time computable
partial functions $\beta, \alpha$ such that 
$ \, f_1 = \beta \circ f_2 \circ \alpha$.
\end{defn} 
Recall polynomial-time {\em many-to-one reduction} that is used for 
languages; it is defined by $L_1 \preccurlyeq_{\sf m:1} L_2$ iff for some 
polynomial-time computable function $\alpha$ and for all $x \in A^*$: 
$x \in L_1$ iff $\alpha(x) \in L_2$.  
This is equivalent to $L_1 = \alpha^{-1}(L_2)$, and also to 
$\chi_{_{L_1}} = \chi_{_{L_2}} \circ \alpha$  (where $\chi_{_{L_j}}$ 
denotes the characteristic function of $L_j$). 
So $L_1 \preccurlyeq_{\sf m:1} L_2$ implies that $\chi_{_{L_1}}$ is simulated
by $\chi_{_{L_2}}$.

Moreover, when we talk about simulations between functions we will always 
use the following 

\medskip

\noindent {\bf Addendum to Definition 5.1.} \ {\it We assume that
$\beta, \alpha \in {\sf fP}$.
For a simulation between two right-ideal morphisms of $A^*$ we assume}
$\beta, \alpha \in {\cal RM}_{|A|}^{\sf P}$.

\bigskip

We can define simulation for monoids in general. For monoids $M_0 \leq M_1$
and $s, t \in M_1$, simulation $s \preccurlyeq t$ is the same thing as  
$s \leq_{{\cal J}(M_0)} t$, i.e., the submonoid ${\cal J}$-order on $M_1$, 
using multipliers in the submonoid $M_0$.

\smallskip

Simulation tells us which functions are harder to compute than 
others, but it does not say anything about the hardness of inverses of
functions. We want a reduction with the property that if a one-way 
function $f_1$ reduces to a function $f_2 \in {\sf fP}$ then $f_2$ is 
also one-way. The intuitive idea is that $f_1$ ``reduces inversively'' to 
$f_2$ iff (1) $f_1$  is simulated by $f_2$, and 
(2) the ``easiest inverses'' of $f_1$ are simulated by the ``easiest 
inverses'' of $f_2$.  But ``easiest inverses'' are difficult to define.
We rigorously define inversive reduction as follows.

\begin{defn} {\bf (inversive reduction).} 
 \ Let $f_1, f_2: A^* \to A^*$ be any partial functions.
We say that $f_1$ {\em reduces inversively} to $f_2$ (notation, 
$f_1 \leqslant_{\sf inv} f_2$) \ iff 

\smallskip

\noindent (1) \ \ $f_1 \preccurlyeq f_2$  \ \ and

\smallskip

\noindent (2) \ \ for every inverse $f_2'$ of $f_2$ there exists an
inverse $f_1'$ of $f_1$ such that $f_1' \preccurlyeq f_2'$ .
\end{defn}
Here, $f_1, f_2, f_1', f_2'$ range over all partial functions $A^* \to A^*$.

\medskip

\noindent The relation $\leqslant_{\sf inv}$ can be defined on monoids 
$M_0 \leq M_1 \leq M_2$ in general: We let $f_1, f_2$ range over $M_1$, and 
let inverses $f_1', f_2'$ range over $M_2$. For simulation $\preccurlyeq$ 
we pick $\leq_{{\cal J}(M_0)}$ (i.e., multipliers are in $M_0$).  
We should assume that $M_1$ is regular within $M_2$ in order to avoid empty
ranges for the quantifiers ``$(\forall f_2')(\exists f_1')$''; otherwise,
when $f_2$ has no inverse in $M_2$, $f_1 \leqslant_{\sf inv} f_2$ is 
trivially equivalent to $f_1 \preccurlyeq f_2$.

\begin{pro} \ $\leqslant_{\sf inv}$ is transitive and reflexive.
\end{pro}
{\bf Proof.} \ Simulation is obviously transitive. Moreover, if
$f_1 \leqslant_{\sf inv} f_2$ and $f_2 \leqslant_{\sf inv} f_3$, then
for each $f_3'$ there exists an inverse
$f_2' = \beta_{23} \circ f_3' \circ \alpha_{23}$, and for $f_2'$ there 
is an inverse $f_1' = \beta_{12} \circ f_2' \circ \alpha_{12}$. 
Then $f_1' = \beta_{12} \circ  \beta_{23} \circ f_3'$ $\circ$
$\alpha_{23} \circ \alpha_{12}$, so $f_3'$ simulates some inverse of $f_1$. 
 \ \ \ $\Box$

\begin{pro} \ If $f_1 \leqslant_{\sf inv} f_2$, 
 \ $f_2 \in {\sf fP}$, and $f_2$ is regular, then $f_1 \in {\sf fP}$ and
$f_1$ is regular.

\noindent Contrapositive: If $f_1, f_2 \in {\sf fP}$, 
$f_1 \leqslant_{\sf inv} f_2$, and $f_1$ is one-way, then $f_2$ is one-way.
\end{pro}
{\bf Proof.} \ The property $f_1 \in {\sf fP}$ follows from simulation. 
If $f_2$ is regular, then it has an inverse $f_2' \in {\sf fP}$, and $f_1$ 
has an inverse $f_1' = \beta \circ f_2' \circ  \alpha$. All the factors 
are in {\sf fP}, so $f_1' \in {\sf fP}$. 
 \ \ \ $\Box$

\begin{defn} \ 
A partial function $g$ is {\em complete} (or {\sf fP}{\em -complete}) with 
respect to $\leqslant_{\sf inv}$ \ iff \ $g \in {\sf fP}$, and for every
$f \in {\sf fP}$ we have $f \leqslant_{\sf inv} g$. In a similar way we can
define ${\cal RM}_2^{\sf P}${\em -complete}.  
\end{defn}
Observation:  If $g$ is {\sf fP}-complete then 
$g \equiv_{\cal J} {\sf id}_{A^*}$.

\begin{pro} \label{ev_complete}
 \ The map ${\sf ev}_{q_2}^C$ is {\sf fP}-{\em complete} with respect to 
inversive reduction.
\end{pro}
{\bf Proof.} \ Any $\phi_w \in {\sf fP}$ with $q_2$-polynomial program $w$ is
simulated by ${\sf ev}_{q_2}^C$; recall formula $(\star)$ in the proof of 
Prop.\ \ref{fPfingen}:

\smallskip

 \ \ \ \ \ $\phi_w \ = \ $
$\pi_{_{2 \, |w'| + 2}}' \circ {\sf contr} \circ {\sf recontr}^{2 \, m}$
$\circ$ ${\sf ev}_{q_2}^C$ $\circ$
${\sf reexpand}^{2 \, m} \circ {\sf expand}$  $\circ$
$\pi_{_{{\sf code}(w) \, 11}}$, 

\smallskip

\noindent where $w' = {\sf co}^{2m+1} \circ {\sf ex}^{2m+1}(w)$.

To prove the inversive property, let ${\sf e}'$ be any inverse of 
${\sf ev}_{q_2}^C$. We apply ${\sf e}'$ to a string of the form \, 
${\sf code}({\sf ex}^{2m+1}(w)) \ 11 \ 0^{N_{2m+1}} \ 11 \ y$,
where $\phi_{{\sf ex}^{2m+1}(w)} \in {\sf fP}^{q_2}$ and 
$y \in {\sf Im}(\phi_w)$.
Thus, $N_{2m+1}$ is at least as large as the time of the computation that 
led to output $y$. Note that we use $2m$ in $N_{2m+1}$ because the input 
that led to output $y$ could be polynomially longer than $y$ (by polynomial 
$q_2$).  Then we have: 

\smallskip

 \ \ \ \ \ 
${\sf e}'\big({\sf code}({\sf ex}^{2m+1}(w))\ 11\ 0^{N_{2m+1}}\ 11\ y \big)$
$ \ = \ $
${\sf code}({\sf ex}^{2m+1}(w)) \ 11 \ 0^{N_{2m+1}} \, 11 \, x_i$, 
 \ \ for some $x_i \in \phi_w^{-1}(y)$.

\smallskip

\noindent We don't care whether and how ${\sf e}'(Z)$ is defined when the
input $Z$ is not of the above form.
Then ${\sf e}'$ simulates an inverse $f'$ of $\phi_w$ defined by 

\smallskip

 \ \ \ \ \ $f'(y) \ = \ $ 
$\pi_{_{2 \, |w'| + 2}}' \circ {\sf contr} \circ {\sf recontr}^{2 \, m}$
$\circ$ ${\sf e}'$ $\circ$
${\sf reexpand}^{2 \, m} \circ {\sf expand}$  $\circ$
$\pi_{_{{\sf code}(w) \, 11}}(y)$

\smallskip

\noindent for all $y \in {\sf Im}(\phi_w)$.  Indeed, $f'(y) = x_i$ 
($\in \phi_w^{-1}(y)$ as above). 

When $y \not\in {\sf Im}(\phi_w)$ the right side of the above formula may
give a value to $f'(y)$; but it does not matter whether and how $f'$ is  
defined outside of ${\sf Im}(\phi_w)$.  Thus, 
every inverse of ${\sf ev}_{q_2}^C$ simulates an inverse of $\phi_w$.
 \ \ \ $\Box$

\medskip

In a similar way one can prove that the generator-based evaluation map 
${\sf ev}_{\Gamma,q}$ (for a large enough polynomial $q$) is complete
in {\sf fP}.

\medskip

\noindent {\bf Notation:} \ For a partial function $f: A^* \to A^*$, the 
{\em set of all inverses} $f': A^* \to A^*$ of $f$ is denoted by 
${\sf Inv}(f)$.

\begin{defn} {\bf (uniform inversive reduction).} 
 \ Let $f, g$ be partial functions.  An inversive reduction 
$f \leqslant_{\sf inv} g$ is called {\em uniform} \ iff
 \ $f \preccurlyeq g$, and $(\exists \beta, \alpha \in {\sf fP})$
$(\forall g' \in {\sf Inv}(g))$ $(\exists f' \in {\sf Inv}(f))$ 
$[ \, f' = \beta \circ g' \circ \alpha \, ]$. So $\beta$ and $\alpha$ only
depend on $f$ and $g$, but not on $g'$ or $f'$.
\end{defn}
We observe that in the proof of Prop.\ \ref{ev_complete} the simulation of
$f'$ by ${\sf e'}$ only depends on $\phi_w$ and ${\sf e}$, but not on 
$f'$ nor on ${\sf e'}$. We conclude:
  
\begin{cor} 
The map ${\sf ev}_{q_2}^C$ is {\sf fP}-complete with respect to 
{\em uniform} inversive reduction.  \ \ \ $\Box$
\end{cor}

Next we study the completeness of the circuit evaluation map 
 ${\sf ev}_{\sf circ}$ (defined at the beginning of Section 4).
Since it is defined in terms of length-preserving circuits, 
${\sf ev}_{\sf circ}$ is itself length-preserving, i.e., it belongs to 
the submonoid of {\em length-preserving} partial functions in
{\sf fP},

\smallskip

 \ \ \  \ \ \ ${\sf fP}_{\sf lp} \ = \ $
$\{ f \in {\sf fP} \ : \ |f(x)| = |x|$ \ for all $x \in {\sf Dom}(f) \}$.

\begin{pro} \label{LevComplLP}
The critical map ${\sf ev}_{\sf circ}$ is complete in the submonoid
${\sf fP}_{\sf lp}$ with respect to inversive reduction.
\end{pro}

\noindent {\bf Proof.} Let $f \in {\sf fP}_{\sf lp}$ be a fixed 
length-preserving partial function, and let $M$ be a fixed deterministic 
polynomial-time Turing machine that computes $f$.

Simulation of $f$ by ${\sf ev}_{\sf circ}$: It is well known that for every 
input length $n$ (of inputs of $f$) one can construct an acyclic circuit 
$C_n$ such that $C_n(x) = f(x)$ for all $x$ of length $n$. The circuit can 
be constructed from $M$ and $n$ in polynomial time (as a function of $n$). 
Let $\alpha(x) = (C_{|x|}, x)$, and let $\beta(C_n, y) = y$, where 
$|y| = n$.
Then $f = \beta \circ {\sf ev}_{\sf circ} \circ \alpha$.
  
Simulation between inverses: Any inverse $e'$ of ${\sf ev}_{\sf circ}$ has 
the form $e'(C,y) = (C, x_i)$ for some $x_i \in C^{-1}(y)$, when
$y \in {\sf Im}(C)$. When $y \not\in {\sf Im}(C)$, $e'(C,y)$ could be any 
pair of bitstrings.
Then an inverse $f'$ of $f$ is obtained by defining 
$f'(y) = \beta \circ e' \circ \alpha(y)$, where $\alpha, \beta$ are as in
the simulation of $f$ (in the first part of this proof).  Indeed, when 
$y \in {\sf Im}(f)$ we have $\alpha: y \longmapsto (C_n, y)$, where $|y| = n$.
Next, $e': (C_n, y) \longmapsto (C_n, x_i)$ for some 
$x_i \in C_n^{-1}(y) = f^{-1}(y)$; recall that we only use length-preserving
circuits, so $|y| = n = |x_i|$.
Finally, $\beta: (C_n, x_i) \longmapsto x_i \in f^{-1}(y)$.  So $f'$ is an 
inverse of $f$ on ${\sf Im}(f)$; outside of ${\sf Im}(f)$, the values of 
$f'$ do not matter.
 \ \ \ $\Box$

\bigskip

\noindent To show completeness of ${\sf ev}_{\sf circ}$ in {\sf fP} (rather
than just in ${\sf fP}_{\sf lp}$), a stronger inversive reduction is needed, 
that overcomes the limitation of length-preservation in ${\sf ev}_{\sf circ}$. 

\bigskip

\noindent {\bf Remark:} 
Circuits are usually generalized to allow the output length to be different 
from the input length. But that would not simplify the completeness proof for
${\sf ev}_{\sf circ}$, 
because the main limitation is that all inputs of a circuit have the same 
length, and all outputs of a circuit have the same length.

\begin{defn} \label{Tsimulation} {\bf (polynomial-time Turing simulation).}
Let $f_1, f_2: A^* \to A^*$ be two partial functions.
By definition, $f_1 \preccurlyeq_{\sf T} f_2$ \ iff \ $f_1$ is computed by a 
deterministic polynomial-time Turing machine that can make oracle calls to 
$f_2$; these can include, in particular, calls on the membership problem of 
${\sf Dom}(f_2)$.
\end{defn}  
In the next proofs we do not need the full power of Turing reductions. 
The following, much weaker reduction, will be sufficient.

\medskip

\noindent {\bf Notation:} Let $L \subseteq A^*$. Then ${\sf fP}^L$ denotes 
the set of all polynomially balanced partial functions computed by 
deterministic polynomial-time Turing machines that can make oracle calls to 
the membership problem of $L$. In particular we will consider 
${\sf fP}^{{\sf Dom}(f)}$ for any partial function $f: A^* \to A^*$. 

\begin{defn} \label{weakTsimulation} {\bf (weak Turing simulation).} 
A {\em weak Turing simulation} of $f_1$ by $f_2$ consists of two partial 
functions $\beta, \alpha$ such that $f_1 = \beta \circ f_2 \circ \alpha$, 
where $\alpha \in {\sf fP}^{{\sf Dom}(f_2)}$ and $\beta \in {\sf fP}$.
The existence of a weak Turing simulation of $f_1$ by $f_2$ is denoted
by $f_1 \preccurlyeq_{\sf wT} f_2$.
\end{defn}
Informally we also write 
$f_1 = \beta \circ f_2 \circ \alpha^{{\sf Dom}(f_2)}$. 
In a weak Turing simulation by $f_2$, only the domain of $f_2$ is repeatedly
queried; $f_2$ itself is called only once, and this call of $f_2$ takes the 
form of an ordinary (not a Turing) simulation.

\begin{defn} \label{inversif} {\bf (inversification of a simulation).}
For any simulation relation $\preccurlyeq_{\sf X}$ between partial functions,
the corresponding {\em inversive reduction} $\leqslant_{\sf inv,X}$ is
defined as follows:

 \ \ $f_1 \leqslant_{\sf inv,X} f_2$ \ \ iff

\smallskip

 \ \   
$f_1 \preccurlyeq_{\sf X} f_2$, and for every inverse $f_2'$ of $f_2$ there 
exists an inverse $f_1'$ of $f_1$ such that $f_1' \preccurlyeq_{\sf X} f_2'$.
\end{defn}
One easily proves:

\begin{pro}
 \ If $\preccurlyeq_{\sf X}$ is transitive then $\leqslant_{\sf inv,X}$ is
transitive.  \ \ \ \ \ \ $\Box$
\end{pro}
Based on this general definition we can introduce {\em polynomial-time 
inversive Turing reductions}, denoted by $\leqslant_{\sf inv,T}$, and 
{\em polynomial-time inversive weak Turing reductions}, denoted by
$\leqslant_{\sf inv,wT}$.
The following is straightforward to prove.

\begin{pro} 
If $f_1 \leqslant_{\sf inv,T} f_2$ then: 

\noindent $\bullet$ \ $f_2 \in {\sf fP}$ implies $f_1 \in {\sf fP}$;

\noindent $\bullet$ \ $f_2 \in {\sf fP}$ and $f_2$ is regular, implies
$f_1$ is regular.
 \ \ \ \ \ $\Box$
\end{pro}
The following shows that $\leqslant_{\sf inv, wT}$ can overcome the 
limitations of length-preservation.

\begin{pro} \label{lpVSall} 
For every $f \in {\sf fP}$ there exists $\ell \in {\sf fP}_{\sf lp}$ 
such that $f \leqslant_{\sf inv,wT} \ell$.
\end{pro}
{\bf Proof.} For any $f \in {\sf fP}$ we define 
$\ell_f \in {\sf fP}_{\sf lp}$ by

\medskip

 \ \ \ \ \ \ \ $\ell_f(0^n 1 \, x) \ = \ \left\{ 
 \begin{array}{ll}
    0^{|x|} 1 \, f(x) \ \ \ \ \ \      & \mbox{if $n = |f(x)|$, } \\
    {\rm undefined}         & \mbox{on all other inputs.}
 \end{array}  \right.
$

\medskip

\noindent Let $p_f(.)$ be a polynomial upper-bound on the time-complexity and 
on the balance of $f$.

\smallskip

\noindent 1. \ Proof that $f \preccurlyeq \ell_f$ (simulation):
 \ We have $f = \beta \circ \ell_f \circ \alpha$, where 
 \ $\alpha(x) = 0^{|f(x)|} 1 \, x$ \ for all $x \in A^*$; and  
 \ $\beta(0^m 1 \, y) = y$ \ for all $y \in A^*$ and all 
 $m \leq p_f(|y|)$ \ ($\beta$ is undefined otherwise).

\smallskip

\noindent  2. \ Proof that for every inverse $\ell'$ of $\ell_f$ there is
an inverse $f'$ of $f$ such that $f' \preccurlyeq_{\sf wT} \ell_f'$ :

\smallskip

Every element of ${\sf Im}(\ell_f)$ has the form $0^m 1 \, y$
where $y \in {\sf Im}(f)$, for some appropriate $m$. More precisely,
$\ell_f^{-1}(0^m \, 1 \, y)$ $ = $ 
$\{ 0^{|y|} \, 1 \, x : x \in f^{-1}(y) \cap A^m\}$.  Hence, 
$0^m \, 1 \, y \in {\sf Im}(\ell_f)$ iff 
$f^{-1}(y) \cap A^m \neq \varnothing$.  Therefore, any inverse $\ell'$ 
satisfies $\ell'(0^m 1 \, y) = 0^{|y|} 1 \, x_i$ for some choice of 
$x_i \in f^{-1}(y) \cap A^m$; we do not care about the values of $\ell'$ 
when its inputs are not in ${\sf Im}(\ell_f)$.
Thus we can define an inverse of $f$ on each $y \in {\sf Im}(f)$ by 

\smallskip

$f'(y) = x_i$, \, for $x_i$ chosen in $f^{-1}(y) \cap A^m$ 

\hspace{.5in} where $m$ is the minimum integer such that 
$f^{-1}(y) \cap A^m \neq \varnothing$. 

\smallskip

\noindent We don't care what $f'(y)$ is when $y \not\in {\sf Im}(f)$.
 
To obtain an inversive weak Turing reduction we need to compute 
$x_i = f'(y)$ from $y$, based on oracle calls to ${\sf Dom}(\ell')$ and
one simulation of $\ell'$.  This computation of $x_i$ is done in two steps:  
 First we compute the minimum $m$ ($= |x_i|$) such that 
$f^{-1}(y) \cap A^m \neq \varnothing$ (see Step 1 below for details). 
 Second, we apply $\ell'$ to compute
$\ell'(0^{|x_i|} 1 \, y) = 0^{|y|} 1 \, x_i$. From this we obtain $x_i$ by 
applying the map $\beta$ defined above (in part 1 of this proof). 
The first step is a Turing reduction to the domain of $\ell'$. 
The second step is a simulation by $\ell'$.  In more detail:

Step 1: 
By input balance we have $|x_i| \leq p_f(|y|)$ when $x_i \in f^{-1}(y)$.
For each $m \in \{0, 1, \ldots, p_f(|y|)\}$, in increasing order, we make 
an oracle call to the membership problem in ${\sf Dom}(\ell')$ with query 
input $0^m 1 \, y$.
If $y \in {\sf Im}(f)$ then the first of these queries with a positive
answer determines $m$, and $0^m 1 \, y$ is returned.

Step 2: To the result $0^m 1 \, y$ of step 1 we apply the functions
$\ell'$ and $\beta$.  This yields $x_i$, which is $f'(y)$. 
Thus, step 2 is just a simulation. 

Togetherm, steps 1 and 2 form a weak polynomial Turing simulation of $f'$ 
by $\ell'$.  
 \ \ \ $\Box$

\begin{cor}
The critical map ${\sf ev}_{\sf circ}$ is {\sf fP}-complete with respect to 
composites of polynomial inversive weak Turing reductions and polynomial
inversive simulation reductions ($\leqslant_{{\sf inv,wT}}$ and 
$\leqslant_{\sf inv}$).
\end{cor}  
{\bf Proof.} For every $f \in {\sf fP}$ we first reduce $f$ to a
length-preserving function $\ell_f$, by Prop.\ \ref{lpVSall}. Then we reduce 
$\ell_f$ to ${\sf ev}_{\sf circ}$ by Prop. \ref {LevComplLP}.
 \ \ \ $\Box$

\bigskip

\noindent {\bf Reduction and completeness in ${\cal RM}_2^{\sf P}$}

\medskip

\noindent The following shows that the encoding that embeds {\sf fP} into
${\cal RM}_2^{\sf P}$ does not make inversion easier.

\begin{pro} \label{code_invred} 
 \ For all $f \in {\sf fP}$ we have $f \leqslant_{\sf inv} f^C$, 
where $\leqslant_{\sf inv}$ is based on simulation in {\sf fP}. 
\end{pro}
{\bf Proof.} Recall the encoding maps
 $(.)_{\#}: x \in \{0,1\}^* \longmapsto x \# \in \{0,1, \#\}^*$, and 
${\sf code}$ which replaces $0, 1, \#$ by respectively $00, 01, 11$, defined
in Section 3; and recall $f^C$ from Def.\ \ref{f_encoding}.  
We now introduce inverses of these maps.
Let ${\sf dec}: {\sf code}(x) \in \{00, 01, 11\}^* $
$\longmapsto x \in \{0,1\}^*$ (undefined outside of $\{00, 01, 11\}^*$), 
and $r: x \# \longmapsto x \in \{0,1\}^*$ (undefined outside of 
$\{0,1\}^*\#$).  Then 

\smallskip

 \ \ \  \ \ \ $f \ = \ $
$r \circ {\sf dec} \circ f^C \circ {\sf code} \circ (.)_{\#}$ .

\smallskip

\noindent Clearly, $(.)_{\#}, \, {\sf code}, \, {\sf dec}, \, r$
$ \in {\sf fP}$.  Hence  $f^C$ simulates $f$. 

For the inversive part of the reduction, let $\varphi'$ be any inverse of
$f^C$; we want to find an inverse $f'$ of $f$ such that 
$f' \preccurlyeq \varphi'$, where $\preccurlyeq$ denotes simulation in 
{\sf fP}. 
Any element of ${\sf Im}(f^C)$ has the form ${\sf code}(s) \, 11 \, t$, 
with $s, t \in \{0,1\}^*$, and $s \in {\sf Im}(f)$. Moreover, if 
${\sf code}(s) \, 11 \, t \in {\sf Im}(f^C)$ then
${\sf code}(s) \, 11 \in {\sf Im}(f^C)$. Let us define $f'$ for any
$s \in {\sf Im}(f)$ by $f'(s) = x_1$ where $x_1$ is such that 
$\varphi'({\sf code}(s) \, 11) = {\sf code}(x_1) \, 11$ 
$ \in (f^C)^{-1}({\sf code}(s) \, 11)$. Then $x_1 \in f^{-1}(s)$. 
In general, finally, we define $f'$ by 

\smallskip 

 \ \ \ \ \ \ $f' \ = \ $
    $r \circ {\sf dec} \circ \varphi' \circ {\sf code} \circ (.)_{\#}$ .

\smallskip

\noindent For $s \in {\sf Im}(f)$ we indeed have then:
$r \circ {\sf dec} \circ \varphi' \circ {\sf code} \circ (.)_{\#}(s) \ = \ $
$r \circ {\sf dec} \circ \varphi'({\sf code}(s) \, 11) \ = \ $
$r \circ {\sf dec}({\sf code}(x_1) \, 11) \ = \ x_1$, \,  where 
$x_1 \in f^{-1}(s)$, as above. So this definition makes $f'$ an inverse of 
$f$ on ${\sf Im}(f)$; hence $f'$ is an inverse of $f$. 
The above formula for $f'$ explicitly shows that $f' \preccurlyeq \varphi'$. 
 \ \ \ $\Box$

\medskip

Let $\equiv_{\sf inv}$ denote $\leqslant_{\sf inv}$-equivalence (i.e.,
$f \equiv_{\sf inv} g$ iff $f \leqslant_{\sf inv} g$ and 
$g \leqslant_{\sf inv} f$).
The $\leqslant_{\sf inv}$-complete functions of {\sf fP} obviously form an
$\equiv_{\sf inv}$-class, and this is the maximum class for the 
$\leqslant_{\sf inv}$-preorder. 
Similarly, the complete functions of ${\cal RM}_2^{\sf P}$ are the maximum 
inversive reducibility class in ${\cal RM}_2^{\sf P}$.
The non-empty regular elements of ${\cal RM}_2^{\sf P}$ also form an 
equivalence class in ${\cal RM}_2^{\sf P}$, and this is the minimum class 
of all non-empty functions, as the following shows:

\begin{pro} \ For every $f, r \in {\cal RM}_2^{\sf P}$ where $r$ is 
regular and $f$ is non-empty, we have $r \leqslant_{\sf inv} f$.
\end{pro}
{\bf Proof.} The simulation $r \preccurlyeq f$ follows from 
${\cal J}^0$-simplicity of ${\cal RM}_2^{\sf P}$.
Let $f'$ be any inverse of $f$ (with $f'$ not necessarily in 
${\cal RM}_2^{\sf P}$). Since $r$ is regular, there is an inverse 
$r' \in {\cal RM}_2^{\sf P}$ of $r$. Since $f'$ is not the empty map there 
exist $x_0, y_0$ with $f'(y_0) = x_0$. Then $(x_0 \leftarrow y_0)$ is 
simulated by $f'$, since
$(x_0 \leftarrow y_0) = {\sf id}_{\{x_0\}} \circ f'$. 
Moreover, $(x_0 \leftarrow y_0)$ is regular and $(x_0 \leftarrow y_0)$ 
belongs to ${\cal RM}_2^{\sf P}$, so $(x_0 \leftarrow y_0)$ simulates $r'$ 
(again by ${\cal J}^0$-simplicity of ${\cal RM}_2^{\sf P}$).  Thus, $f'$ 
simulates $r'$.
 \ \ \ $\Box$

\begin{pro} \ In both {\sf fP} and ${\cal RM}_2^{\sf P}$:  The 
$\equiv_{\cal D}$-relation is a refinement of $\equiv_{\sf inv}$.
\end{pro}
{\bf Proof.} Is suffices to prove that both $\equiv_{\cal R}$ and 
$\equiv_{\cal L}$ refine $\equiv_{\sf inv}$.
We will prove that $f \equiv_{\cal R} g$ implies $f \equiv_{\sf inv} g$.
(The same reasoning works for $\equiv_{\cal L}$.) 
When $f \equiv_{\cal R} g$, there exist $\alpha, \beta \in {\sf fP}$ (or 
$\in {\cal RM}_2^{\sf P}$) such that $f = g \, \alpha$ and $g = f \, \beta$. 
So, $f$ and $g$ simulate each other.

For any inverse $f'$ of $f$ we have $f = f \, f' \, f $
$ = g \, \alpha f' f$. Right-multiplying by $\beta$ we obtain
$g = g \, \alpha f' \, g$, hence $\alpha f'$ is an inverse of $g$, and 
$\alpha f'$ is obviously simulated by $f'$. So, $g$ inversely reduces to 
$f$. Similarly, $f$ inversely reduces to $g$.
 \ \ \ $\Box$

\section{\bf The polynomial hierarchy}

\noindent The classical polynomial hierarchy for languages is defined by
$\Sigma_1^{\sf P} = {\sf NP}$, $\Pi_1^{\sf P} = {\sf coNP}$, and 
for all $k > 0$: $\Sigma_{k+1}^{\sf P} = {\sf NP}^{\Sigma_k^{\sf P}}$
(i.e., all languages accepted by nondeterministic Turing machines with oracle 
in $\Sigma_k^{\sf P}$, or equivalently, with oracle in $\Pi_k^{\sf P}$); and
$\Pi_{k+1}^{\sf P} = ({\sf coNP})^{\Sigma_k^{\sf P}}$ 
$\big( = {\sf co}({\sf NP}^{\Sigma_k^{\sf P}})\big)$. Moreover,
${\sf PH} \ = \ \bigcup_k \Sigma_k^{\sf P}$ \ \ ($\subseteq {\sf PSpace}$).

\medskip

\noindent {\sf Polynomial hierarchy for functions:}

\smallskip

${\sf fP}^{\Sigma_k^{\sf P}}$ \, consists of all polynomially balanced 
partial functions $A^* \to A^*$ computed by deterministic polynomial-time 
Turing machines with oracle in $\Sigma_k^{\sf P}$ (or equivalently, with 
oracle in $\Pi_k^{\sf P}$);

\smallskip

${\sf fP}^{\sf PH}$ \, consists of all polynomially balanced partial 
functions $A^* \to A^*$ computed by deterministic polynomial-time Turing 
machines with oracle in {\sf PH}. Equivalently, 
${\sf fP}^{\sf PH} = \, \bigcup_k \, {\sf fP}^{\Sigma_k^{\sf P}}$.

\smallskip

Moreover, ${\sf fPSpace}$ \, consists of all polynomially balanced
partial functions (on $A^*$) computed by deterministic polynomial-space 
Turing machines.

\smallskip

We can also define a polynomial hierarchy over ${\cal RM}_2^{\sf P}$.

\begin{pro} \label{fPNP_inv} \ Every $f \in {\sf fP}$ has an inverse in 
${\sf fP}^{\Sigma_1^{\sf P}}$, and
every $f \in {\sf fP}^{\Sigma_k^{\sf P}}$ has an inverse in
${\sf fP}^{\Sigma_{k+1}^{\sf P}}$. The monoids 
${\sf fP}^{\sf PH}$ and ${\sf fPSpace}$ are regular.
\end{pro}
{\bf Proof.} \ The following is an inverse of $f$:

\smallskip

 \ \ \ \ \ $f'_{\sf min}(y) \ = \ \left\{
   \begin{array}{ll}
    {\sf min}(f^{-1}(y))    & \mbox{if $y \in {\sf Im}(f)$, } \\
    {\rm undefined}         & \mbox{otherwise,} 
   \end{array}  \right.
$

\smallskip

\noindent where ${\sf min}(S)$ denotes the minimum of a set of strings $S$ 
in dictionary order (or alternatively in length-lexicographic order).
To show that $f'_{\sf min} \in {\sf fP}^{\sf NP}$ when $f \in {\sf fP}$ we 
first observe that for any fixed $f \in {\sf fP}$ the following problems are 
in {\sf NP}: 

\noindent (1) On input $y \in A^*$, decide whether $y \in {\sf Im}(f)$.

\noindent (2) Fix $u \in A^*$; on input  $y \in A^*$, decide whether 
$y \in f(u \, A^*)$ (i.e., decide whether there exists $x \in u \, A^*$ such 
that $f(x) = y$).

When $y \not\in {\sf Im}(f)$ then it doesn't matter what value we choose 
for $f'_{\sf min}(y)$; we choose $f'_{\sf min}(y)$ to be undefined then.

Here is an ${\sf fP}^{\sf NP}$-algorithm for computing $f'_{\sf min}(y)$.
It is a form of binary search in the sorted list $A^*$, that ends when 
a string in $f^{-1}(y)$ has been found. A growing prefix $z$ of 
$x = f'_{\sf min}(y)$ is constructed; at each step we query whether $z$ can 
be extended by a 0 or a 1; i.e., we ask whether $y \in f(z0 A^*)$; we don't 
need to ask whether $y \in f(z1 A^*)$ too, since we tested already that 
$y \in {\sf Im}(f)$.  Oracle calls are denoted by angular brackets 
$\langle \ldots \rangle$, and $\varepsilon$ denotes the empty word. 

\smallskip

\noindent Algorithm for $f'_{\sf min}$ on input $y:$

\smallskip

{\tt if $\langle y \in {\sf Im}(f) \rangle$ then

 \ \ \ $z := \varepsilon$;
 
 \ \ \ while $\langle z \not\in f^{-1}(y) \rangle$ do \hspace{1in}
    // assume $y \in f(z A^*)$

 \ \ \  \ \ \ if $\langle y \in f(z0 A^*) \rangle$ then 
    $z := z 0$;   

 \ \ \  \ \ \ else $z := z 1$; 

 \ \ \ output $z$.
}   %%% \tt

\smallskip

\noindent One can prove in a similar way that when 
$f \in {\sf fP}^{\Sigma_k^{\sf P}}$ then $f'_{\sf min} \in$
${\sf fP}^{\Sigma_{k+1}^{\sf P}}$: In that case the problems (1) and (2)
above are in ${\sf NP}^{ \Sigma_k^{\sf P}} = \Sigma_{k+1}^{\sf P}$.   

The regularity of ${\sf fP}^{\sf PH}$ follows immediately from the 
fact about ${\sf fP}^{\Sigma_k^{\sf P}}$ for each $k$.
The regularity of ${\sf fPSpace}$ holds because the above algorithm 
can be carried out in ${\sf fPSpace}$.  
 \ \ \ $\Box$

\bigskip

The above algorithm is similar to the proofs in the literature  that 
{\sf P} $\neq$ {\sf NP} iff one-way functions exist; see e.g.\
\cite{HemaOgi} p.\ 33.

\medskip

In the proof of Prop.\ \ref{fPNP_inv} we used the function $f'_{\sf min}$.
In a similar way, by using ${\sf max}(f^{-1}(y))$ one can define
$f'_{\sf max} \in {\sf fP}^{\Sigma_1^{\sf P}}$, which is also an inverse of 
$f$. Yet more inverses can be defined:  for any positive integer $i$ let

\smallskip

 \ \ \ \ \ $ f'_i(y) \ = \ \left\{
\begin{array}{ll}
i^{\rm th} \ {\rm word \ in} \ f^{-1}(y) \ \ \ 
                        & \mbox{if $y \in {\sf Im}(f)$,} \\
{\rm undefined}         & \mbox{otherwise.}
\end{array}  \right.
$

\smallskip

\noindent Here, ``$i^{\rm th}$ word'' refers to the dictionary order; also, 
when $i > |f^{-1}(y)|$, the $i^{\rm th}$ word is taken to be the maximum 
word in $f^{-1}(y)$. Then $f'_i$ is an inverse of $f$ and 
$f'_i \in {\sf fP}^{\Sigma_1^{\sf P}}$; note that $i$ is fixed for each 
function $f'_i$.

\begin{pro} \ For any {\em {\sf fP}-critical} partial function 
$f \in {\sf fP}$ we have:
 \ $f$ is one-way \, iff \, $f'_{\sf min} \not\in {\sf fP}$.
Similarly,  \,  $f$ is one-way \, iff \, $f'_{\sf max} \not\in {\sf fP}$
\, iff \, $(\exists i > 0)[ \, f'_i \not\in {\sf fP} \, ]$.
\end{pro}
{\bf Proof.}  Since $f'_{\sf min}$ is an inverse of $f$, the direction 
``$\Rightarrow$'' is obvious by the definition of one-way function. 
Conversely, we saw that if $f \in {\sf fP}$ then 
$f'_{\sf min} \in {\sf fP}^{\Sigma_1^{\sf P}}$. 
If $f'_{\sf min} \not\in {\sf fP}$ then 
${\sf fP} \neq {\sf fP}^{\Sigma_1^{\sf P}}$, hence ${\sf P} \neq {\sf NP}$, 
hence one-way functions exist.
Then any {\sf fP}-critical function $f$ is one-way.
 \ \ \ $\Box$

\medskip

\noindent Recall that a partial function $f$ is called {\sf fP}-critical iff
$f \in {\sf fP}$ and the existence of one-way functions implies that $f$ is
one-way. In particular, {\sf fP}-complete functions (with respect to inversive
reduction) are {\sf fP}-critical.
An interesting consequence of the above Proposition is that now we do not 
only have critical functions, but these functions also have {\em critical 
inverses}.

\begin{defn} \ Let $f$ be an {\sf fP}-critical function.  We say that an 
inverse $f'$ of $f$ is a {\em critical inverse} of $f$ \, iff \,    
$f' \not\in {\sf fP}$ implies that $f$ is one-way.
\end{defn}

\begin{cor}
For the {\sf fP}-critical function ${\sf ev}_{\sf circ}$, the inverses 
$({\sf ev}_{\sf circ})'_{\sf min}$, $({\sf ev}_{\sf circ})'_{\sf max}$
and $({\sf ev}_{\sf circ})'_i$ are critical inverses.
 \ \ \ $\Box$ 
\end{cor}
Thus, to decide whether ${\sf P} \neq {\sf NP}$ it suffices to consider 
one function, and one of its inverses.

\begin{pro} \label{polyHierFinGen}  \ For each $k \geq 1$ the monoid 
${\sf fP}^{\Sigma_k^{\sf P}}$ is finitely generated, but not finitely 
presented. The monoid ${\sf fPSpace}$ is also finitely generated, 
but not finitely presented.

The monoid ${\sf fP}^{\sf PH}$ is not finitely generated, unless the
polynomial hierarchy collapses.
\end{pro}
{\bf Proof.} The proof for ${\sf fPSpace}$ is similar to the proof that we 
gave for {\sf fP} in Prop.\ \ref{notFinPres}. 

For ${\sf fP}^{\Sigma_k^{\sf P}}$, let $Q_k$ be any 
$\Sigma_k^{\sf P}$-complete problem; we can assume that all oracle calls
are calls to $Q_k$. Then every $f \in {\sf fP}^{\Sigma_k^{\sf P}}$ has a
program which is like an {\sf fP}-program, but with oracle calls to $Q_k$ 
added. 
For every polynomial $q \geq q_2$, an evaluation function ${\sf ev}_q^{Q_k}$ 
for ${\sf fP}^{\Sigma_k^{\sf P}}$ can then be designed; in the computation of 
${\sf ev}_q^{Q_k}(w, x)$, oracle calls to $Q_k$ are made whenever the program 
$w$ being executed makes calls to $Q_k$. Then, 
${\sf ev}_q^{Q_k}(w, x) = (w, \phi_w(x))$.
By using ${\sf ev}_q^{Q_k}$ the proof for ${\sf fP}^{\Sigma_k^{\sf P}}$
is similar to the proof of Prop.\ \ref{notFinPres}.  

If ${\sf fP}^{\sf PH}$ were finitely generated then let $m$ be the lowest
level in the hierarchy that contains a finite generating set. Then 
${\sf fP}^{\sf PH} \subseteq {\sf fP}^{\Sigma_m^{\sf P}}$.
 \ \ \ $\Box$

\bigskip

Instead of using all of ${\sf fP}^{\sf NP}$ to obtain inverses for the 
elements of {\sf fP}, we could simply adjoin inverses to {\sf fP} (within
${\sf fP}^{\sf NP}$). 
It turns out that it suffices to adjoin just one inverse 
$e' \in {\sf fP}^{\sf NP}$ of a function $e$ that is {\sf fP}-complete for 
$\leqslant_{\sf inv}$.

\medskip

\noindent 
{\bf Notation:} For a semigroup $S$ and a subset $W \subseteq S$, the 
subsemigroup of $S$ generated by $W$ is denoted by $\langle W \rangle_S$.
For any $h \in {\sf fP}^{\sf NP}$, we denote 
$\langle {\sf fP} \cup \{h\} \rangle_{{\sf fP}^{\sf NP}}$ by 
${\sf fP}[h]$ (called {\em ``{\sf fP} with $h$ adjoined''}).
So, \ \   
${\sf fP} \ \subseteq \ {\sf fP}[h] \ \subseteq \ {\sf fP}^{\sf NP}$.

\begin{pro}
Let $g \in {\sf fP}$ be any function that is {\sf fP}-complete with respect
to $\leqslant_{\sf inv}$, and let $g'$ be any inverse of $g$ such that 
$g' \in {\sf fP}^{\sf NP}$.
Then the subsemigroup ${\sf fP}[g']$ of ${\sf fP}^{\sf NP}$ contains at least
one inverse of each element of {\sf fP}. 
\end{pro}
{\bf Proof.} From the assumption that $g$ is complete it follows that

\smallskip

 \ \ \   \ \ \   \ \ \   $(\forall f \in {\sf fP})$
$(\forall g' \in {\sf Inv}(g) \cap {\sf fP}^{\sf NP})$
$(\exists f' \in {\sf Inv}(f))$ $(\exists \beta, \alpha \in {\sf fP})$ 
$[ \, f' = \beta \, g' \, \alpha \, ]$. 

\smallskip

\noindent So for any fixed $g' \in {\sf Inv}(g) \cap {\sf fP}^{\sf NP}$, 
every $f \in{\sf fP}$ has an inverse of the form 
$f' = \beta \, g' \, \alpha$, for some $\beta, \alpha \in {\sf fP}$ (that
depend on $f'$).  Hence $f' \in {\sf fP}[g']$.  
\ \ \ $\Box$

\bigskip

\noindent {\bf Observations:} 

\smallskip

\noindent 1.
We saw in the proof of Prop.\ \ref{polyHierFinGen} that ${\sf fP}^{\sf NP}$ 
has complete elements with respect to simulation. For any 
${\sf fP}^{\sf NP}$-complete element $h$ we have
${\sf fP}^{\sf NP} = {\sf fP}[h]$. This raises the question: 
Is ${\sf fP}^{\sf NP} \neq {\sf fP}[g']$, when $g' \in {\sf fP}^{\sf NP}$ and
$g'$ is an inverse of an element $g$ that is {\sf fP}-complete 
(for $\leqslant_{\sf inv}$)? 
In one direction we have:

\smallskip

{\it If there exists $g$ which is {\sf fP}-complete with respect 
to $\leqslant_{\sf inv}$, and an inverse $g' \in {\sf fP}^{\sf NP}$ such that 
${\sf fP}^{\sf NP} \neq {\sf fP}[g']$, then ${\sf P} \neq {\sf NP}$.
}

Indeed, if ${\sf fP}^{\sf NP} \neq {\sf fP}[g']$ then
${\sf fP} \subseteq {\sf fP}[g'] \varsubsetneq {\sf fP}^{\sf NP}$, hence
${\sf fP} \neq {\sf fP}^{\sf NP}$, hence ${\sf P} \neq {\sf NP}$.

\medskip

\noindent 2. {\it If there exist $g_1, g_2$ (not necessarily different) that
are {\sf fP}-complete with respect to $\leqslant_{\sf inv}$, and inverses
$g'_1, g'_2 \in {\sf fP}^{\sf NP}$ of $g_1$, respectively $g_2$, such that
${\sf fP}[g'_1] \neq {\sf fP}[g'_2]$, then ${\sf P} \neq {\sf NP}$.
}

Indeed, by contraposition, if ${\sf P} = {\sf NP}$ then 
${\sf fP} = {\sf fP}^{\sf NP}$, hence $g'_1, g'_2 \in {\sf fP}$. Then
${\sf fP}[g'_1] = {\sf fP} = {\sf fP}[g'_2]$.

\medskip

\noindent 3. {\it The following two statements are equivalent: 
 \ (1) \ ${\sf P} \neq {\sf NP}$; 
 \ (2) \ there exist $g$ which is {\sf fP}-complete with respect to
$\leqslant_{\sf inv}$, and an inverse $g' \in {\sf fP}^{\sf NP}$ such that
${\sf fP} \neq {\sf fP}[g']$.
}

Indeed, if such a $g$ and $g'$ exist then $g$ is a one-way function,
hence ${\sf P} \neq {\sf NP}$.
If for such a $g$ and $g'$ we have ${\sf fP} = {\sf fP}[g']$, then $g$ is
an {\sf fP}-complete function which is not one-way, hence one-way functions
do not exist.

\bigskip

\noindent {\bf Other monoids:} 

\smallskip

\noindent (1) We have: \, ${\sf fP}^{\sf PSpace} = {\sf fPSpace}$. 

Indeed, the monoid ${\sf fP}^{\sf PSpace}$ consists of polynomially 
balanced functions that are polynomial-time computable, with calls to 
{\sf PSpace} oracles. 
The monoid {\sf fPSpace} consists of polynomially balanced functions that 
are polynomial-space computable (hence they might use exponential time).
Obviously, ${\sf fP}^{\sf PSpace} \subseteq {\sf fPSpace}$. But the 
converse holds too, since the polynomially many output bits of a function
in {\sf fPSpace} can be found one by one, by a polynomial number of calls
to {\sf PSpace} oracles.

\medskip

\noindent (2) We define {\sf fLog} (``functions in log-space'') to 
consist of the polynomially balanced partial functions that are computable 
in deterministic log space.  {\sf fLog} is closed under composition (see 
\cite{HU}), and ${\sf fLog} \subseteq {\sf fP}$.

If ${\sf P} \neq {\sf NP}$ then {\sf fLog} is non-regular;
more strongly, in that case {\sf fLog} contains one-way functions (with 
no inverse in {\sf fP}). Indeed, the 3{\sc cnf} 
satisfiability problem reduces to the inversion of the map
$(B, \alpha) \mapsto (B, B(\alpha))$, where $B$ is any boolean formula in 
3{\sc cnf}, and $\alpha$ is a truth-value assignment for $B$. It is not
difficult to prove that this map is in {\sf fLog} when $B$ is in 3{\sc cnf}. 
One of the referees observed that {\sf fLog} is regular iff {\sf NP} $=$
{\sf L}, i.e., the class of languages accepted in deterministic log-space. 

\medskip

\noindent (3) We define {\sf fLin} (``functions in linear time'') to 
consist of the linearly balanced partial functions that are computable 
in deterministic linear time. {\sf fLin} is closed under composition, and
it is non-regular \ iff \ ${\sf P} \neq {\sf NP}$. More strongly, if ${\sf P}
\neq {\sf NP}$ then {\sf fLin} contains one-way functions (with no inverse 
in {\sf fP}); this is proved by padding arguments.

\bigskip

\bigskip

\noindent {\bf Acknowledgement:}  This paper benefitted greatly form 
corrections offered by the referees.

%%%%%%%%%%%%%%%%%%%%%%%%%%%%%%%%%%%%%%%%%%%%%%%%%%%%%%%%%%%%%%

%% \small

%%%%%%

\vspace{.5in}

J.C.\ Birget 

Dept.\ of Computer Science

Rutgers University {\bf --} Camden 

Camden, New Jersey 

\smallskip

{\tt birget@camden.rutgers.edu}

\end{document}